\documentclass{article}
\usepackage{amsfonts}
\usepackage{latexsym}
\usepackage{amssymb}
\usepackage{amsmath}
\usepackage[mathscr]{euscript}
\usepackage{mathrsfs}
\usepackage{marvosym}
\usepackage{fancyhdr}
 \usepackage{cite}
\usepackage[numbers,sort&compress]{natbib}
\usepackage{color}
\newtheorem{lemma}{Lemma}[section]
\newtheorem{corollary}{Corollary}[section]
\newtheorem{proposition}{Proposition}[section]
\newtheorem{theorem}{Theorem}[section]
\newtheorem{remark}{Remark}[section]
\newtheorem{definition}{Definition}[section]
\newtheorem{example}{Example}[section]
 
\def\cvd{\hfill $\Box$}

\title{{{\bf\large  Existence  of Weak Pareto Efficient Solutions of a Vector Optimization Problem under a Closed Constraint Set }}\thanks{This
work was supported by the National Natural Science Foundation of China (11871059), the Applied Basic Project of Sichuan Province (2020YJ0111), the Innovation Team of Department of Education of Sichuan Province (16TD0019), the Meritocracy Research Funds of China West Normal University (17YC379) and the Fundamental Research Funds of China West Normal University (493006).}}

\author {\small Danyang Liu\thanks{Key Laboratory of Optimization Theory and Applications at China West Normal University of Sichuan Province, School of Mathematics and Information, China West Normal University, 637009 Nanchong, Sichuan, P. R. China (Corresponding author, e-mail: dyliu@cwnu.edu.cn)},\, Jun Li\thanks{Key Laboratory of Optimization Theory and Applications at China West Normal University of Sichuan Province, School of Mathematics and Information, China West Normal University, 637009 Nanchong, Sichuan, P. R. China (junli@cwnu.edu.cn)}\, and\, Giandomenico Mastroeni \thanks{Department of Computer Science, University of  Pisa, Largo B. Pontecorvo 3, 56127 Pisa,  Italy (gmastroeni@di.unipi.it)}
}

\begin{document}
\date{}
\maketitle{}

{\bf Abstract.}  In this paper, we investigate the nonemptiness of weak Pareto efficient solution set for a  class of nonsmooth vector optimization problems on a nonempty closed constraint set without any boundedness and convexity assumptions. First, we obtain a new property  concerning the nonemptiness of weak Pareto efficient solution sets for these vector optimization problems. Then, under the condition of  weak section-boundedness from below, we establish relationships between the notions of properness, Palais-Smale, weak Palais-Smale and M-tameness conditions with respect to some index sets for the restriction of the vector mapping on the constraint set. Finally, we present some new necessary and sufficient conditions  for the nonemptiness of the weak Pareto efficient solution set for a general class of nonsmooth vector optimization problems.

{\bf Key words.} Existence  of weak solutions,  M-Tame, Properness, Palais-Smale conditions.

{\bf AMS subject classifications.} 90C29, 90C30, 49J30

\def\cvdeq{\hfill \sqcap \hskip-5.5pt \sqcup}

\def\vs{\vskip.3truecm}

\def\inte{{\rm \ int\ }}
\def\cl{{\rm \ cl\ }}

\section{Introduction}
\setcounter{equation}{0}

Throughout the paper,  $\mathbf{R}^{n}$ denotes the $n$-dimensional Euclidean space with  the norm $\|\cdot\|$ and the inner  product $\langle \cdot,\cdot\rangle$,  and $\mathbf{R}^{n}_{+}:=\{x=(x_1,\cdots,x_n)\in \mathbf{R}^{n}:x_i \ge 0, i=1, \cdots, n\}$. The closed unit ball in $\mathbf{R}^n$ is denoted by $\mathbf{B}^n$ and the natural numbers set is denoted by $\mathbf{N}$.
We are interested in the following vector optimization problem:
$$ \text{VOP}(K, f): \mbox{Min}_{x\in K}f(x),$$
where $f:=(f_{1},  \dots, f_{s}): \mathbf{R}^{n} \to \mathbf{R}^{s}$ is  a locally Lipschitz mapping, that is, for any $x\in \mathbf{R}^n$, there exist a real number $L>0$ and a neighborhood $U(x)$ of $x$ such that  $\|f_i(x'')-f_i(x')\|\leq L \|x''-x'\|$  for any $x',x''\in U(x)$, and the constraint set $K\subseteq \mathbf{R}^{n}$ is a nonempty closed set (not necessarily a convex  or semi-algebraic set  \cite{HHV,BR}). Let  $I:=\{s_1,s_2,\dots,s_q\}$, with $s_1<s_2<\cdots<s_q$, be an arbitrary nonempty index  subset of the set $\{1,2,\dots,s\}$ and assume $|I|=q$ stands for the number of the elements in the set $I$. Let $f_I: \mathbf{R}^n\to \mathbf{R}^q$ be given by $f_I(x):=(f_{s_1}(x),f_{s_2}(x),\dots,f_{s_q}(x))$ for any $x\in\mathbf{R}^n$.

\par Recall that a point $x^{*}\in K$ is a Pareto efficient solution of $\text{VOP}(K, f)$, if $$f(x)-f(x^{*})\notin - \mathbf{R}^{s}_{+}\backslash \{0\}, \quad \forall x\in K,$$
and   $x^{*}\in K$ is a weak Pareto efficient solution of VOP$(K, f)$, if $$f(x)-f(x^{*})\notin - \mbox{ int }\mathbf{R}^{s}_{+}, \quad \forall x\in K,$$
where $\mbox{ int } \mathbf{R}^{s}_{+}:=\{x=(x_1,\cdots,x_s)\in \mathbf{R}^{s}:x_i > 0, i=1, \cdots, s\}$. The Pareto efficient solution set and the weak Pareto efficient solution set of $\text{VOP}(K, f)$ are denoted by ${\rm SOL}^{S}(K, f)$ and ${\rm SOL}^{W}(K, f)$, respectively. Clearly, ${\rm SOL}^{S}(K, f)\subseteq {\rm SOL}^{W}(K, f)$.  When $s=1$,  $\text{VOP}(K, f)$ collapses to a scalar optimization problem on $K$, denoted by $\mbox{SOP}(K, f)$.

As ${\rm SOL}^{S}(K, f)$  may be empty, in this paper, we are concerned with the fundamental problems   revolved around the existence of weak Pareto efficient solutions of $\text{VOP}(K, f)$. Some related results on the existence of Pareto efficient solutions for vector optimization problems were investigated under some restrictive assumptions. Indeed, in \cite{Kim2}, when $K=\mathbf{R}^n$ and $f$ is a polynomial mapping, by powerful methods of semi-algebraic geometry and polynomial optimization theory, Kim, Pham and Tuyen studied the existence of Pareto efficient solutions under the properness, the Palasi-Smale, the weak Palais-Smale and the M-tame definitions. When the constraint set $K$ is a  convex semi-algebraic  set and  $f$ is a convex polynomial mapping, Lee, Sisarat and Jiao \cite{LGJ} proved that ${\rm SOL}^{S}(K, f)$  is nonempty if and only if  the image $f(K)$ of $f$ has a nonempty bounded section. When $K$ is a  closed  semi-algebraic set (without convexity) and  $f$ is a polynomial mapping, Duan et al. \cite{Duan} obtained the existence of Pareto efficient solutions for polynomial vector optimization problems by employing  semi-algebraic theory under the section-boundedness from below condition. In \cite{LDY2}, Liu, Huang and Hu studied the existence of Pareto efficient solutions for a polynomial vector optimization problem when the weak section-boundedness from below, the convenience and non-degeneracy conditions were satisfied. When $K$ is a  closed  set and  $f$ is a polynomial mapping, Liu, Hu and Fang \cite{LDY} studied the solvability of a regular polynomial vector optimization problem by using  asymptotic analysis tools.  By means of  variational analysis, Kim et al. \cite{Kim3} investigated the nonemptiness of Pareto efficient solution set for a general class of nonsmooth vector optimization problems on a nonempty closed constraint set without any boundedness assumption.

 Some related results on existence of weak Pareto efficient solutions were also studied under some restrictive assumptions. In the case of unconstrained vector optimization problems, existence theorems for weak Pareto efficient solutions and  Pareto efficient solutions  have been obtained in \cite{BT1,BT2,MK2} by using appropriate set-valued extensions of the Ekeland variational principle under the boundedness from below and Palais-Smale condition assumptions. When the constraint set is convex, Flores-Baz\'{a}n and Vera \cite{Flores2} obtained  an existence  theory for the convex multiobjective optimization problems by employing an asymptotic approach.  In \cite{Flores}, Flores-Baz\'{a}n established existence of the weak efficient solutions for finite dimensional VOP allowing the solution set to be unbounded. In \cite{LDY2}, Liu, Huang and Hu proved that a  convex polynomial vector optimization problem  admits a weak Pareto efficient solution if and only if the weak section-boundedness from below condition is satisfied. Based on a new coercivity notion for vector-valued functions, Guti\'{e}rrez and L\`{o}pez \cite{Guti} studied vector optimization problems with solid non-polyhedral convex ordering cones by using the Gerstewitz scalarization function, asymptotic analysis and a regularization of the objective function without any convexity or quasiconvexity assumption. For  other results related to these topics see also  the books \cite{Jahn,Luc,Sawaragi}, the papers \cite{Borwein,Deng1,Deng2,Huy,Lara,Luc2} and references therein.

 \par In this paper, we consider the nonemptiness of weak Pareto efficient solution set for a general class of  vector optimization problems. 	
Our main contributions are the following:
 \begin{itemize}
  	\item We consider  nonsmooth  vector optimization problems with an arbitrary nonempty closed constraint set, without any polynomial, convexity and boundedness assumptions.
 	\item  We obtain some new properties on the nonemptyness of weak Pareto efficient solution set of these  vector optimization problems,  proving new  necessary and sufficient conditions for the existence of weak Pareto efficient solutions.
	
	\item Our approach is mainly based on tools of variational analysis and generalized differentiation, but not on semi-algebraic theorem.
 \end{itemize}

The  paper is organized as follows. In Section 2, we present some notations and preliminary results.  In Section 3, we obtain some new  characterizations of the nonemptyness of weak Pareto efficient solution set of $\text{VOP}(K, f)$. In Section 4, we establish relationships between the notions of properness, Palais-Smale, weak Palais-Smale and M-tameness conditions with respect to some index set for the restriction of the vector mapping on the constraint set. In Section 5, we obtain some new necessary and sufficient conditions  for the nonemptiness of weak Pareto efficient solution set for a  general class of  nonsmooth vector optimization problems and finally, 	 Section 6 is devoted to concluding remarks.

\section{Preliminaries}
\setcounter{equation}{0}
In this section, we recall the main concepts and results that will be used in the sequel.

\subsection{Definitions of section-boundedness from below}

First, we recall the definition of  {\it a section of a set}. Let $A$ be a subset of $\mathbf{R}^s$ and $\bar t\in \mathbf{R}^s$. The set $A\bigcap (\bar t -\mathbf{R}^s_+)$ is called a section of $A$ at $\bar t$ and denoted by $\lbrack A \rbrack_{\bar t}$. The section $\lbrack A \rbrack_{\bar t}$ is said to be bounded if there exists $a\in \mathbf{R}^s$ such that $\lbrack A \rbrack_{\bar t}\subseteq a +\mathbf{R}^s_+$, see, e.g. \cite{Kim2,Jahn,0810ha,Borwein}. Based on this concept, the following definitions are introduced.

\begin{definition}\label{strong section2}
Let $A$ be a nonempty subset of $\mathbf{R}^s$ and $\bar t\in \mathbf{R}^s$. The section $\lbrack A \rbrack_{\bar t}$ is said to be bounded with respect to $I$ at $\bar t$, if there  exists a nonempty index set $I\subseteq \{1,2,\dots,s\}$, with $q:=|I|$, such that
	the set $\lbrack A \rbrack^I_{\bar t}:=\{(x_{s_1},x_{s_2},\dots,x_{s_q})\in\mathbf{R}^q: (x_1,x_2,\dots,x_s)\in \lbrack A \rbrack_{\bar t}\}$, is bounded,  i.e.,
	$\lbrack A \rbrack^I_{\bar t}\subseteq \bar a +\mathbf{R}^q_+$, for some $\bar a\in \mathbf{R}^q$.
\end{definition}
\begin{remark}\label{Rem:1}
	By definition, if $I=\{1,2,\dots,s\}$, then the boundedness with respect to $I$ at $\bar t$ of the section $\lbrack A \rbrack_{\bar t}$ is equivalent to the boundedness  of the section $\lbrack A \rbrack_{\bar t}$.
\end{remark}

 Next, we recall the definition of {\it weak section-boundedness from below} for a vector-valued mapping, that is, a vector-valued mapping $f:=(f_{1},  \dots, f_{s}): \mathbf{R}^{n} \to \mathbf{R}^{s}$ is said to be weakly section-bounded from below on $K$, if there exists $a\in \mathbf{R}^{s}$ such that $f(x)-a\notin -\mbox{ int }\mathbf{R}^{s}_{+}$ for any $x\in K_{\bar{x}}$, where $K_{\bar{x}}:=\{x\in K: f_{i}(x)\leq f_{i}(\bar{x}), i=1, 2, \cdots, s\}$, see e.g. \cite{LDY2}. Then by Definition 3 and Proposition 1 in \cite{LDY2}, we  know that $f$ is weakly section-bounded from below on $K$ if and only if there exist $x^*\in K$ and $i_{0}\in \{1, 2, \cdots, s\}$ such that $f_{i_{0}}$ is bounded from below on $K_{x^*}$. Motivated by this, we  introduce the following definition.

\begin{definition}\label{weak section-boundedness}
	Let $f: \mathbf{R}^{n}\rightarrow \mathbf{R}^{s}$ be a vector-valued mapping. Then $f$ is said to be weakly section-bounded from below on $K$ with respect to $I\subseteq \{1,2,\dots,s\}$ and $x_0\in K$, if the section $\lbrack f(K) \rbrack_{f(x_0)}$ is bounded with respect to $I$ at $f(x_0)$, that is, there exist a nonempty index set $I\subseteq \{1,2,\dots,s\}$ and $x_0\in K$ such that the set $\lbrack f(K) \rbrack^I_{f(x_0)}$ is bounded, where $\lbrack f(K) \rbrack^I_{f(x_0)}:=\{y_I\in \mathbf{R}^{|I|}|(y_1,y_2,\dots,y_s)\in \lbrack f(K) \rbrack_{f(x_0)}\}$.
\end{definition}

\begin{remark}\label{0805weak}

 Since $[f(K)]^I_{f(x_0)}=\{f_I(x)\in \mathbf{R}^{|I|}|(f_1(x),f_2(x),\dots,f_s(x))\in \lbrack f(K) \rbrack_{f(x_0)}\}=\{f_I(x)\in \mathbf{R}^{|I|}|x\in K, f_i(x)\leq f_i(x_0), i=1,2,\dots,s\}=\{f_I(x)\in \mathbf{R}^{|I|}|x\in K_{x_0}, f_i(x)\leq f_i(x_0), i\in I\}=[f_I(K_{x_0})]_{f_I(x_0)}$, then $f$ is  weakly section-bounded from below on $K$ with respect to $I\subseteq \{1,2,\dots,s\}$ and $x_0\in K$, if and only if the section $\lbrack f_I(K_{x_0}) \rbrack_{f_I(x_0)}$ is bounded.

  It is known that $f$ is said to be section-bounded from below on $K$, if the section $\lbrack f(K) \rbrack_{f(x_0)}$ is bounded (see e.g. \cite[Definition 2.2]{LDY2}). By Definition \ref{weak section-boundedness}, if $I=\{1,2,\dots,s\}$, then the weak section-boundedness from below implies the section-boundedness from below.  Moreover,  we observe   that  weak section-boundedness from below on $K$ of $f$ is a necessary condition  for the nonemptiness of ${\rm SOL}^{W}(K, f)$ (see \cite[Remark 2.7]{LDY2}).
\end{remark}

\subsection{Normal cones and subdifferentials}

We recall the notions of  normal cone of a closed set and of subdifferential of a real-valued function, more details can be found in \cite{MK1,Rock2}.

\begin{definition}
	Let $K$ be a subset of $\mathbf{R}^n$ and  $\bar x$ in $K$.
	\begin{itemize}
		\item[\rm(i)]  The regular normal cone $\hat{N}(\bar x; K)$ to $K$ at $\bar x$ consists of all vectors $v\in \mathbf{R}^n$ satisfying
	    $$\langle v, x-\bar x \rangle\leq o(\|x-\bar x\|), \quad \forall x\in K,$$
	    where the term $o(\|x-\bar x\|)$ satisfies $\frac{o(\|x-\bar x\|)}{\|x-\bar x\|}\to 0$ as $x\to \bar x$ with $x\in K\backslash\{\bar x\}$.
	   \item[\rm(ii)]  The limiting normal cone $N(\bar x; K)$ to $K$ at $\bar x$ consists of all vectors $v\in \mathbf{R}^n$ such that there are sequences $x_k\to \bar x$ with $x_k\in K$ and $v_k\to v$ with $v_k\in \hat{N}(x_k; K)$ as $k\to +\infty$.
	\end{itemize}
\end{definition}

\begin{definition}
	Consider a function $\varphi: \mathbf{R}^n\to \mathbf{R}$ and a point $\bar x\in \mathbf{R}^n$. The  subdifferential of $\varphi$ at $\bar x$ is defined by
	$$\partial\varphi(\bar x):=\{v\in\mathbf{R}^n\vert (v, -1)\in N((\bar x, \varphi(\bar x)); \mbox{epi }\varphi) \},$$
	where $N((\bar x, \varphi(\bar x)); \mbox{epi }\varphi)$ is  the normal cone to the epigraph $\mbox{epi }\varphi$ of $\varphi$, where
$$\mbox{epi }\varphi:=\{(x, y)\in \mathbf{R}^n\times\mathbf{R}\vert \varphi(x)\leq y\}.$$
\end{definition}

 We also need the following results for establishing relationships between the properness, Palais-Smale conditions, and M-Tameness of $f$ on $K$ with respect to some index set and the necessary and the sufficient condition for the nonemptiness of $\text{SOL}^W(K,f)$. First, we recall  classical subdifferential formulas and optimality conditions of convex  and variational analysis.

\begin{lemma}\cite{Rock1} \label{classical subdifferential}
For any $x^*\in \mathbf{R}^n$, we have that the subdifferential $\partial g$ of the function $g: \mathbf{R}^n\to \mathbf{R}, x\mapsto\|\cdot-x^*\|$ is as follows:
\begin{itemize}
	\item[\rm(i)] If $x\neq x^*$, then $\partial(\|\cdot-x^*\|)(x)=\frac{x-x^*}{\|x-x^*\|}$.
	\item[\rm(ii)] If $x= x^*$, then $\partial(\|\cdot-x^*\|)(x)=\mathbf{B}^n$.
\end{itemize}
\end{lemma}


\begin{lemma}\label{optimality condition}\cite[Corollary 6.6]{MK2}
	Let $\varphi_i: \mathbf{R^n\to \mathbf{R}} (i=0,1,\dots, s+r)$ be locally Lipschitz functions around $x^*\in K$ and let the set $K\subseteq\mathbf{R}^n$ be locally closed around $x^*\in K$. If $x^*$ is a local minimizer of the function $\varphi_0$ on the set
	$$\{x\in K| \varphi_i(x)\leq 0,i=1,2,\dots,s,\varphi_i(x)=0, i=s+1,s+2,\dots,s+r\},$$
	then there exists a nonzero collection of multipliers $\{\alpha_0,\alpha_1,\dots,\alpha_{s+r}\}\in \mathbf{R}^{s+r+1}$ with $\alpha_i\geq 0, i=0,1,\dots,s$, such that
	$$0\in \sum_{i=0}^{s}\alpha_i\partial \varphi_i(x^*)+\sum_{i=s+1}^{s+r}\alpha_i\partial^0 \varphi_i(x^*)+N(x^*,K),$$
	$$\alpha_i\varphi_i(x^*)=0, i=1,2,\dots,s,$$
	where $\partial^0 \varphi_i(x^*)$ is the symmetic subdifferential of $\varphi_i$ at $x^*$ defined by
	$$\partial^0 \varphi_i(x^*):=-\partial(-\varphi_i)(x^*)\bigcup \partial\varphi_i(x^*),i=s+1,s+2,\dots,s+r.$$
\end{lemma}

It is known that, if $\varphi_i$ is strictly differentiable at $x^*$, then $\partial^0 \varphi_i(x^*)=\partial\varphi_i(x^*)=\{\nabla \varphi_i(x^*)\}$.

\begin{lemma}\label{sum fule}\cite[Theorem 3.36]{MK1}
	Let $\varphi_i: \mathbf{R^n\to \mathbf{R}} (i=1,2,\dots, s)$ be locally Lipschitz functions around $x^*\in \mathbf{R}^n$. Then we have the subdifferential sum rule
	$$\partial(\sum_{i=1}^{s}\varphi_i)(x^*)\subseteq \partial \varphi_1(x^*)+\partial \varphi_2(x^*)+\cdots+\partial \varphi_s(x^*).$$
\end{lemma}

\begin{lemma}\label{maximum rule}\cite[Theorem 3.46]{MK1}
	Let $\varphi_i: \mathbf{R^n\to \mathbf{R}} (i=1,2,\dots, s)$ be locally Lipschitz functions around $x^*\in \mathbf{R}^n$. Then the maximum function $\varphi(x):=\max_{1\leq i \leq s}\varphi_i(x) (x\in \mathbf{R}^n)$
	is locally Lipschitz  around $x^*\in \mathbf{R}^n$ and the following inclusion holds
	$$\partial\varphi(x^*)\subseteq\{\sum_{i\in I(x^*)}\alpha_i\partial \varphi_i(x^*)| \alpha_i\geq 0, \sum_{i\in I(x^*)}\alpha_i=1\},$$
	where  $I(x^*):=\{i\in \{1,2,\dots,s\}| \varphi(x^*)=\varphi_i(x^*)\}$.
\end{lemma}

\section{New characterizations of the nonemptiness of weak Pareto efficient solutions set for $\text{VOP}(K,f)$}
\setcounter{equation}{0}

In this section, we  introduce some new results  concerning the existence of weak Pareto efficient solutions for $\text{VOP}(K,f)$, which will be used in proving a characterization of the nonemptiness of  $\text{SOL}^W(K,f)$. Let $f=(f_{1},  \dots, f_{s}): \mathbf{R}^{n} \to \mathbf{R}^{s}$ be a vector-valued mapping (not necessarily a locally Lipschitz mapping) and $K$ be a nonempty subset of $\mathbf{R}^n$ (not necessarily a closed set). The following result holds.

\begin{lemma}\label{2023801}
	Let $f$ be a vector-valued mapping as above. Then for any $x'\in K$, ${\rm SOL}^{W}(K_{x'}, f)\subseteq {\rm SOL}^{W}(K, f)$ and ${\rm SOL}^{S}(K_{x'}, f)\subseteq {\rm SOL}^{S}(K, f)$,  where $K_{x'}:=\{x\in K| f_i(x)\leq f_i(x'), i=1,2,\dots, s\}$.
\end{lemma}
{\it Proof}\,\, It suffices to prove the first conclusion. Similarly, we can show the second one.
For any $x'\in K$, let $x''\in {\rm SOL}^{W}(K_{x'}, f)$. For any $x\in K$, if  $x\in K_{x'}$, then it follows from $x''\in {\rm SOL}^{W}(K_{x'}, f)$ that $f(x)-f(x'')\notin - \mbox{ int }\mathbf{R}^{s}_{+}$; if $x\notin K_{x'}$, then there exists $i'\in \{1,2,\dots,s\}$ such that
 $f_{i'}(x)-f_{i'}(x') > 0$. Since $x''\in {\rm SOL}^{W}(K_{x'}, f)$,  $f_{i}(x'')\leq f_{i}(x')$ for any $i=1,2,\dots, s$. It follows that $f_{i'}(x)-f_{i'}(x'') > 0$ and so $f(x)-f(x'')\notin - \mbox{ int }\mathbf{R}^{s}_{+}$. As a consequence, $f(x)-f(x'')\notin - \mbox{ int }\mathbf{R}^{s}_{+}$ for any $x\in K$, that is, $x''\in {\rm SOL}^{W}(K, f)$. This completes the proof. \cvd

Now, we introduce new  results concerning the existence of weak Pareto efficient solutions for $\text{VOP}(K,f)$.

\begin{proposition}\label{property of nonempty}
	Assume that ${\rm SOL}^{W}(K, f)\neq\emptyset$. Then there exist $i_0\in \{1,\dots,s\}$ and $x_0\in K$ such that $f_{i_0}(x)\equiv f_{i_0}(x_0)$ for any $x\in K_{x_0}$, where $K_{x_0}:=\{x\in K| f_i(x)\leq f_i(x_0), i=1,2,\dots,s\}$. Moreover, $x_0\in {\rm SOL}^{W}(K_{x_0}, f)\subseteq {\rm SOL}^{W}(K, f)$.
\end{proposition}
{\it Proof}\,\, Since ${\rm SOL}^{W}(K, f)\neq\emptyset$, let $\bar x\in {\rm SOL}^{W}(K, f)$. Consider
$$K_{\bar x}:=\{x\in K| f_i(x)\leq f_i(\bar x), i=1,2,\dots, s\}.$$
Then $K_{\bar x}\neq\emptyset$. If $f_{1}(x)\equiv f_{1}(\bar x)$ for any $x\in K_{\bar x}$, then we only let $i_0:=1$. It is easy to see that $\bar x\in {\rm SOL}^{W}(K_{\bar x}, f)$. By Lemma \ref{2023801}, we have ${\rm SOL}^{W}(K_{\bar x}, f)\subseteq {\rm SOL}^{W}(K, f)$. It follows that $\bar x\in {\rm SOL}^{W}(K, f)$. If not, then there exists $x_1\in K_{\bar x}$ such that $f_{1}(x_1)< f_{1}(\bar x)$. Consider
$$K_{x_1}:=\{x\in K| f_i(x)\leq f_i(x_1), i=1,2,\dots, s\}.$$
Then $K_{x_1}\neq\emptyset$. Since $x_1\in K_{\bar x}$, we have $K_{x_1}\subseteq K_{\bar x}$. If $f_{2}(x)\equiv f_{2}(x_1)$ for any $x\in K_{x_1}$, then let $i_0:=2$ and so $x_1\in {\rm SOL}^{W}(K_{x_1}, f)$. By Lemma \ref{2023801}, we have ${\rm SOL}^{W}(K_{x_1}, f)\subseteq {\rm SOL}^{W}(K, f)$, and so $x_1\in {\rm SOL}^{W}(K, f)$.  If not, then there exists $x_2\in K_{x_1}$ such that $f_{2}(x_2)< f_{2}(x_1)$. Consider
$$K_{x_2}:=\{x\in K| f_i(x)\leq f_i(x_2), i=1,2,\dots, s\}.$$
Then $K_{x_2}\neq\emptyset$. Since $x_2\in K_{x_1}$, we have $K_{x_2}\subseteq K_{x_1} \subseteq K_{\bar x}$. If $f_{3}(x)\equiv f_{3}(x_2)$ for any $x\in K_{x_2}$, then let $i_0:=3$ and so $x_2\in {\rm SOL}^{W}(K_{x_2}, f)$. By Lemma \ref{2023801}, we have ${\rm SOL}^{W}(K_{x_2}, f)\subseteq {\rm SOL}^{W}(K, f)$, and so $x_2\in {\rm SOL}^{W}(K, f)$.  If not, then there exists $x_3\in K_{x_2}$ such that $f_{3}(x_3)< f_{3}(x_2)$. Consider the set
$$K_{x_3}:=\{x\in K| f_i(x)\leq f_i(x_3), i=1,2,\dots, s\}.$$
Then $K_{x_3}\neq\emptyset$. Since $x_3\in K_{x_2}$, we have $K_{x_3}\subseteq K_{x_2} \subseteq K_{x_1} \subseteq K_{\bar x}$.
\par Repeating this process,  we can know that if $f_{s}(x)\equiv f_{s}(x_{s-1})$ for any $x\in K_{x_{s-1}}$, then let $i_0:=s$ and so $x_{s-1}\in {\rm SOL}^{W}(K_{x_{s-1}}, f)$. By Lemma \ref{2023801}, we have ${\rm SOL}^{W}(K_{x_{s-1}}, f)\subseteq {\rm SOL}^{W}(K, f)$, and so $x_{s-1}\in {\rm SOL}^{W}(K, f)$.  If not, then there exists $x_s\in K_{x_{s-1}}$ such that $f_{s}(x_s)< f_{s}(x_{s-1})$. Consider  $K_{x_s}:=\{x\in K| f_i(x)\leq f_i(x_s), i=1,2,\dots, s\}$. Then $K_{x_s}\neq\emptyset$. Since $x_s\in K_{x_{s-1}}$, we have $K_{x_s}\subseteq K_{x_{s-1}} \subseteq\cdots\subseteq K_{x_1} \subseteq K_{\bar x}$. By the definition of the set $K_{x_i},i=1,2,\dots,s$, we can obtain the following relationships:
	\begin{equation*}
		\begin{split}
		f_s(x_s)<f_s(x_{s-1})\leq f_s(x_{s-2})\leq f_s(x&_{s-3})\leq\cdots \leq f_s(x_1)\leq f_s(\bar x),\\
		f_{s-1}(x_s)\leq f_{s-1}(x_{s-1})< f_{s-1}(x_{s-2})&\leq f_{s-1}(x_{s-3})\leq \cdots \leq f_{s-1}(x_1)\leq f_{s-1}(\bar x),\\
		f_{s-2}(x_s)\leq f_{s-2}(x_{s-1})\leq f_{s-2}(x_{s-2})&< f_{s-2}(x_{s-3})\leq \cdots \leq f_{s-2}(x_1)\leq f_{s-2}(\bar x),\\
		&\vdots\\
		f_3(x_s)\leq f_3(x_{s-1})\leq  \cdots \leq f_3(x_3) < &f_3(x_2)\leq f_3(x_1)\leq f_3(\bar x),\\
		f_2(x_s)\leq f_2(x_{s-1})\leq  \cdots \leq f_2(x_3)\leq &f_2(x_2)<f_2(x_1)\leq f_2(\bar x),\\
		f_1(x_s)\leq f_1(x_{s-1})\leq  \cdots \leq f_1(x_3)\leq &f_1(x_2)\leq f_1(x_1)<f_1(\bar x).
		\end{split}
	\end{equation*}
Then, $f_i(x_s)<f_i(\bar x)$ for any $i\in \{1,2,\dots,s\}$. It follows from $x_s\in K_{x_s}\subseteq K_{x_{s-1}} \subseteq\cdots\subseteq K_{x_1} \subseteq K_{\bar x}$ that $\bar x\notin {\rm SOL}^{W}(K, f)$, which  contradicts that $\bar x\in {\rm SOL}^{W}(K, f)$. Therefore, there exist $i_0:=s$ and $x_0:=x_{s-1}$ in $K$ such that $f_{i_0}(x)\equiv f_{i_0}(x_0)$ for any $x\in K_{x_0}$. Since $x_{s-1}\in {\rm SOL}^{W}(K_{x_{s-1}}, f)$ and ${\rm SOL}^{W}(K_{x_{s-1}}, f)\subseteq {\rm SOL}^{W}(K, f)$, we have $x_{s-1}\in {\rm SOL}^{W}(K, f)$. This completes the proof.
\cvd
\begin{remark}\label{I nonempty}
	By Proposition \ref{property of nonempty}, we can easily see that ${\rm SOL}^{W}(K, f)\neq\emptyset$ if and only if there exist $i_0\in \{1,2,\dots,s\}$ and $x_0\in K$ such that $f_{i_0}(x)\equiv f_{i_0}(x_0)$ for any $x\in K_{x_0}$. We consider the following index set
	\begin{equation}\label{L1}
	I:=\{i\in \{1,2,\dots,s \}| f_{i}(x)\equiv f_{i}(x_0),\forall x\in K_{x_0}\}.
	\end{equation}
	Then ${\rm SOL}^{W}(K, f)\neq\emptyset$ if and only if there exists $x_0\in K$ such that $I\neq\emptyset$. In particular, if $I=\{1,2,\dots,s\}$, then we can get a stronger result as follows  (see e.g. \cite{Kim3}).
\end{remark}

\begin{proposition}\label{I=123s}
	$x_0\in{\rm SOL}^{S}(K, f)\neq\emptyset$ if and only if  $x_0\in K$ and $I=\{1,2,\dots,s\}$, where $I$ is defined by \eqref{L1}.
\end{proposition}
{\it Proof}\,\,
``$\Rightarrow$''.\,\,Let $x_0\in {\rm SOL}^{S}(K, f)$. Consider the set $K_{x_0}:=\{x\in K| f_i(x)\leq f_i(x_0), i=1,2,\dots,s\}$. Clearly, for all $i\in \{1,\dots,s\}$, we have $f_{i}(x)\equiv f_{i}(x_0)$ for any $x\in K_{x_0}$. Thus $I=\{1,2,\dots,s\}$.

``$\Leftarrow$''.\,\,Let $x_0\in K$ be such that $I=\{1,2,\dots,s\}$. Then ${x_0\in {\rm SOL}^{S}(K_{x_0}, f)}$ and from Lemma \ref{2023801} it follows that $x_0\in {\rm SOL}^{S}(K, f)\neq\emptyset$. This completes the proof. \cvd

\begin{remark}
	From \cite{Kim3}, we  know that Proposition \ref{I=123s} plays an important role in proving a characterization of the nonemptiness of Pareto efficient solution set for a general class of nonsmooth  vector optimization problems.
\end{remark}

By Remark \ref{I nonempty} we have the following result.
\begin{proposition}\label{P1}
The following assertions are equivalent:
    \begin{itemize}
	\item[\rm(i)]  ${\rm SOL}^{W}(K, f)\neq\emptyset$.
	\item[\rm(ii)] There exists $x_0\in K$ such that the index set $I=\{s_1,s_2,\cdots, s_q\}$ as defined by \eqref{L1}, is nonempty and $x_0\in {\rm SOL}^{S}(K, f_I)$, where $f_I:=(f_{s_1},f_{s_2},\dots,f_{s_q})$.
\end{itemize}
\end{proposition}
{\it Proof}\,\, Assume that (i) holds. Then we have observed in Remark \ref{I nonempty} that there exists $x_0\in K$ such that $I\ne \emptyset$. We now prove  that $x_0\in {\rm SOL}^{S}(K_{x_0}, f_I)$, i.e., the system
$$f_i(x)-f_i(x_0) \le 0, i\in I  \ {\rm and} \ f_j(x)-f_j(x_0) < 0  \ {\rm for \ some} \ j\in I, i\ne j, \quad x\in K_{x_0},$$
is  impossible. Indeed, by the definition of the set $I$ and $K_{x_0}$,  if $\bar x$ is a solution of the previous system, then $\bar x\in K_{x_0}$ with $f_j(\bar x)-f_j(x_0) < 0,$ \ for \ some $ j\in I$,
which contradicts that for $j\in I, x\in K_{x_0}$, $f_j(x)=f_j(x_0)$.
By Lemma \ref{2023801}, we obtain that $x_0\in {\rm SOL}^{S}(K, f_I)$.

Vice-versa, if (ii) holds, then  by definition of weak Pareto efficient solution, we can easily prove that $x_0\in {\rm SOL}^{W}(K, f)$.
In fact, if $x_0\in {\rm SOL}^{S}(K, f_I)$, then  there is no $ x\in K$ such that
$$f_i( x)-f_i(x_0) \le 0, \ i\in I \ {\rm and}\   f_j(x)-f_j(x_0) < 0  \ {\rm for \ some} \ j\in I, i\ne j, $$
which implies that the system
\begin{equation}\label{S1}
f_i( x)-f_i(x_0) < 0, \ i\in \{1,2, ...,s\},\  x\in K,
\end{equation}
is impossible, i.e., (i) holds. This completes the proof.
 \cvd

By the proof of Proposition \ref{P1} we also obtain the following result.

\begin{corollary}\label{C1}
The following assertions are equivalent:
    \begin{itemize}
	\item[\rm(i)]  ${\rm SOL}^{W}(K, f)\neq\emptyset$.
	\item[\rm(ii)] There exists $x_0\in K$ such that the index set $I$ as in \eqref{L1} is nonempty and $x_0\in {\rm SOL}^{S}(K_{x_0}, f_I)$, where $f_I:=(f_{s_1},f_{s_2},\dots,f_{s_q})$ and $I:=\{s_1,s_2,\cdots, s_q\}\subseteq \{1,2,\cdots,s\}$.
	\end{itemize}
\end{corollary}

{\bf Proof.} The implication (i) $\Rightarrow$ (ii) is proved in the previous proposition.

Vice-versa, if (ii) holds, then by $x_0\in {\rm SOL}^{S}(K_{x_0}, f_I)$ and Lemma \ref{2023801}, we obtain that $x_0\in {\rm SOL}^{S}(K, f_I)$. Therefore, by Proposition \ref{P1}, we have that (i) holds. \cvd

\section{Properness, Palais-Smale conditions, and M-Tameness with respect to  an index set}
\setcounter{equation}{0}

In this section, let $f:=(f_1,f_2,\dots,f_s): \mathbf{R}^{n}\rightarrow \mathbf{R}^{s}$ be a locally Lipschitz mapping, let $K$ be a nonempty closed subset in $\mathbf{R}^{n}$ and let  $I\subseteq\{1,2,\dots,s\}$ a nonempty index set.  In order to prove the closedness of a section of $f(K)$, which is difficult to verify, we consider  the notion of properness of the restricted mapping with respect to some index set $I$. Then, under the weak section-boundedness from below assumption, we shall discuss the relationships between the properness, Palais-Smale conditions, and M-Tameness of $f$ on $K$ with respect to the index set $I$. In addition, we obtain the closedness of a bounded section of $f(K)$.

First, we introduce the new notion  of properness with respect to the  index set $I$ for the restricted mapping $f|_K$ of $f$ on $K$.

\begin{definition}\label{Iproper}
Let  $I\subseteq\{1,2,\dots,s\}$ be a nonempty index set. We say that
	\begin{itemize}
		\item[\rm(i)] The restricted mapping $f|_K$ of $f$ on $K$ is said to be proper with respect to $I$ at sublevel $y_0\in \mathbf{R}^s$, if
	$$\forall \{x_k\}\subseteq K, \|x_k\|\to +\infty, f(x_k)\leq y_0\Longrightarrow \|f_I(x_k)\|\to +\infty \mbox{ as } k\to +\infty.$$
	   \item[\rm(ii)]  The restricted mapping $f|_K$ is said to be proper  with respect to $I$, if it is proper with respect to $I$ at every sublevel $y_0\in \mathbf{R}^s$.
	\end{itemize}
\end{definition}

\begin{remark}
	 Similarly to Remark 3.1 in \cite{Kim3}, when $s=1$ and $f$ is bounded from below, the properness with respect to the  index set $I$ of the restricted mapping $f|_K$ is equivalent to the coercivity of $f|_K$. When $s\geq 2$, we  know that the properness with respect to the some index set $I$ of the restricted mapping $f|_K$ is weaker than $\mathbf{R}^s_+$-zero-coercivity of $f$ on $K$ (see e.g. \cite[Definition 3.1]{Kim3}).
\end{remark}

\begin{remark}\label{08051}

 Given the vector function $f:=(f_1,...,f_s)$, let $f_I=(f_i)_{i\in I}$, the vector function with components $f_i, i\in I$. We observe that  $f$ is
proper with respect to $I$ at sublevel $f(x_0)\in \mathbf{R}^s$, if and only if $ f_I$ is proper on $K_{x_0}$ at $f_I(x_0)$ in the sense of \cite[Definition 3.2]{Kim3}.
	In particular, when $I=\{1,2,\dots,s\}$, then Definition \ref{Iproper} reduces to \cite[Definition 3.2]{Kim3}.
	
\end{remark}


\begin{definition}\label{IKT}
	  For any nonempty index set $I\subseteq \{1,2,\dots,s\}$ and $y_0\in (\mathbf{R}\cup\{\infty\})^s$, define the following sets:
$$\widetilde{K}^{I}_{\infty,\leq y_0}(f, K):=\{y\in\mathbf{R}^{|I|}\vert \exists \{x_k\}\subseteq K, f(x_k)\leq y_0, \|x_k\|\to +\infty, f_{I}(x_k) \to y \mbox{ and } \nu(x_k)\to 0 \mbox{ as } k\to +\infty\},$$
$$K^{I}_{\infty,\leq y_0}(f, K):=\{y\in\mathbf{R}^{|I|}\vert \exists \{x_k\}\subseteq K, f(x_k)\leq y_0, \|x_k\|\to +\infty, f_{I}(x_k) \to y  \mbox{ and } \|x_k\|\nu(x_k)\to 0 \mbox{ as } k\to +\infty\},$$
$$\mbox{and} \ T^{I}_{\infty,\leq y_0}(f, K):=\{y\in\mathbf{R}^{|I|}\vert \exists \{x_k\}\subseteq \Gamma(f, K), f(x_k)\leq y_0, \|x_k\|\to +\infty \mbox{ and } f_I(x_k) \to y \mbox{ as } k\to +\infty\},$$
where $\nu:\mathbf{R}^n\to \mathbf{R}\cup\{+\infty\}$ is the extended Rabier function defined by
$$\nu(x):=\inf\{\|\sum_{i=1}^{s}\alpha_i v_i+\omega\|\vert v_i\in \partial f_i(x),i\in\{1,2,\dots,s\}, \omega\in N(x;K), \alpha=(\alpha_1,\alpha_2,\dots,\alpha_s)\in \mathbf{R}^s_+, \sum_{i=1}^{s}\alpha_i=1\},$$
and the tangency variety of $f$ on $K$ defined by
\begin{equation*}
\begin{split}
	&\Gamma(f, K):=\{x\in K\vert v_i\in \partial f_i(x),i\in\{1,2,\dots,s\}, \exists (\alpha, \mu)\in \mathbf{R}^s_+\times \mathbf{R} \mbox{ with }  \sum_{i=1}^{s}\alpha_i+|\mu|=1 \mbox{ such that } \\
		&0\in  \sum_{i=1}^{s}\alpha_i v_i+\mu x+N(x; K)\}.
\end{split}
		\end{equation*}
\end{definition}

\begin{remark}
	Let $I\subseteq \{1,2,\dots,s\}$ and $y_0\in (\mathbf{R}\cup\{\infty\})^s$. Clearly, by the definitions, the inclusion  $K^I_{\infty, \leq y_0}(f, K)\subseteq \widetilde{K}^I_{\infty, \leq y_0}(f, K)$ holds. By \cite[Proposition 3.2]{Kim5}, when $f$ is a polynomial mapping and $K=\mathbf{R}^n$, we have that the inclusion  $T^I_{\infty, \leq y_0}(f, K)\subseteq K^I_{\infty, \leq y_0}(f, K)$ holds. By \cite[Theorem 4.1]{Duan}, when $f$ is a polynomial mapping and $K$ is a closed semi-algebraic set  satisfying regularity at infinity, then the inclusion  $T^I_{\infty, \leq y_0}(f, K)\subseteq K^I_{\infty, \leq y_0}(f, K)$ also holds. However, it is worth noting that if $f$ is not polynomial, then the inclusion  $T^I_{\infty, \leq y_0}(f, K)\subseteq K^I_{\infty, \leq y_0}(f, K)$ may not hold: see, e.g., \cite[Example 3.1]{Kim3}, where $f:\mathbf{R}\to \mathbf{R}$ is defined by $f(x):=\sin x$, $\forall x\in \mathbf{R}$, $K=\mathbf{R}$ and  $y_0=0$. We can see that $0\in T^I_{\infty, \leq y_0}(f, K)\backslash K^I_{\infty, \leq y_0}(f, K)$; thus, we have that $T^I_{\infty, \leq y_0}(f, K)\not\subseteq K^I_{\infty, \leq y_0}(f, K)$.
\end{remark}

\begin{remark}
	In particular, if $I=\{1,2,\dots,s\}$, then $\widetilde{K}^{I}_{\infty,\leq y_0}(f, K)$, $K^{I}_{\infty,\leq y_0}(f, K)$, and $\mbox{and} \ T^{I}_{\infty,\leq y_0}(f, K)$
	reduce to the following sets (see, e.g., \cite{Kim3}):
	$$\widetilde{K}_{\infty,\leq y_0}(f, K):=\{y\in\mathbf{R}^s\vert \exists \{x_k\}\subseteq K, f(x_k)\leq y_0, \|x_k\|\to +\infty, f(x_k) \to y, \mbox{ and } \nu(x_k)\to 0 \mbox{ as } k\to +\infty\},$$
$$K_{\infty,\leq y_0}(f, K):=\{y\in\mathbf{R}^s\vert \exists \{x_k\}\subseteq K, f(x_k)\leq y_0, \|x_k\|\to +\infty, f(x_k) \to y, \mbox{ and } \|x_k\|\nu(x_k)\to 0 \mbox{ as } k\to +\infty\},$$
$$T_{\infty,\leq y_0}(f, K):=\{y\in\mathbf{R}^s\vert \exists \{x_k\}\subseteq \Gamma(f, K), f(x_k)\leq y_0, \|x_k\|\to +\infty, \mbox{ and }f (x_k) \to y \mbox{ as } k\to +\infty\},$$
\end{remark}

\vskip0.2cm

\begin{definition}\label{IPPss}
	Let $I\subseteq \{1,2,\dots,s\}$ be a nonempty index set and $y_0\in (\mathbf{R}\cup\{\infty\})^s$. We say that
	\begin{itemize}
		\item[\rm(i)] $f|_{K}$ satisfies the Palais-Smale condition with respect to $I$ at the sublevel $y_0$ if
		$$\widetilde{K}^{I}_{\infty,\leq y_0}(f, K)=\emptyset.$$
		\item[\rm(ii)] $f|_{K}$ satisfies the weak Palais-Smale condition with respect to $I$ at the sublevel $y_0$ if
		$$K^{I}_{\infty,\leq y_0}(f, K)=\emptyset.$$
		\item[\rm(iii)] $f|_{K}$ M-tame with respect to $I$ at the sublevel $y_0$ if
		$$T^{I}_{\infty,\leq y_0}(f, K)=\emptyset.$$
	\end{itemize}
\end{definition}

\begin{remark}\label{0811PS}
	In particular, if $I=\{1,2,\dots,s\}$, then (i)-(iii) of Definition \ref{IPPss} reduce to the Palais-Smale,  weak Palais-Smale  and M-tame condition, that is, $\widetilde{K}_{\infty,\leq y_0}(f, K)=\emptyset$, $K_{\infty,\leq y_0}(f, K)=\emptyset$ and $T_{\infty,\leq y_0}(f, K)=\emptyset$, (see  \cite[Definition 3.3]{Kim3}).
\end{remark}
	
 From the definitions, it follows  that the properness of the restricted mapping $f|_K$ of $f$ on $K$ with respect to $I$ at sublevel $y_0\in \mathbf{R}^s$ yields
	$$\widetilde{K}^{I}_{\infty,\leq y_0}(f, K)=K^{I}_{\infty,\leq y_0}(f, K)=T^{I}_{\infty,\leq y_0}(f, K)=\emptyset.$$
The converse does not hold in general, see, e.g. \cite{Kim3}.



In what follows, we give necessary and sufficient conditions for the properness with respect to $I$ at the sublevel $y_0$ under the weak section-boundedness from below condition, by slightly modifying \cite[Theorem 3.1]{Kim3}.
\begin{theorem}\label{relationship}
    Assume that $f$ is  weakly section-bounded from below on $K$ with respect to the nonempty index set $I\subseteq\{1,2,\dots,s\}$ and $x_0\in K$. Then the following assertions are equivalent:
    \begin{itemize}
	\item[\rm(i)] $f|_{K}$ is proper with respect to $I$ at the sublevel $f(x_0)$.
	\item[\rm(ii)] $f|_{K}$ satisfies the Palais-Smale condition with respect to $I$ at the sublevel $f(x_0)$.
	\item[\rm(iii)] $f|_{K}$ satisfies the weak Palais-Smale condition with respect to $I$ at the sublevel $f(x_0)$.
  \item[\rm(iv)] $f|_{K}$ satisfies M-tame condition with respect to $I$ at the sublevel $f(x_0)$.
    \end{itemize}
Moreover, the set $\{x\in K|f(x)\leq f(x_0)\}$ and the section $\lbrack f(K) \rbrack_{f(x_0)}$ are compact if any of the conditions (i)-(iv) is fulfilled.
\end{theorem}
{\it Proof}\,\,  By the definition of weak section-boundedness from below, there exist a nonempty index set $I\subseteq \{1,2,\dots,s\}$ and $x_0\in K$ such that the set $\lbrack f(K) \rbrack^I_{f(x_0)}$ is bounded. Note that ${\rm (i)}\Rightarrow {\rm(ii)}, {\rm(ii)}\Rightarrow{\rm(iii)}$ and ${\rm(i)}\Rightarrow {\rm(iv)}$ are obvious.

Now, we first prove ${\rm(iii)}\Rightarrow {\rm(i)}$. Suppose, on the contrary, that $f|_{K}$ is not proper with respect to $I$ at the sublevel $f(x_0)$. Let $i_0\in I$. Consider the following set:
\begin{equation}\label{0802Kx}
	K^{i_0}_{x_0}:=\{x\in K|f_i(x)\leq f_i(x_0),i=1,2,\dots,i_0-1,i_0+1,\dots,s\}.
\end{equation}
Then it is obvious that the set $K^{i_0}_{x_0}$ is nonempty and unbounded. Let
\begin{equation}\label{0802b}
	b:=\lim_{x\in K^{i_0}_{x_0},\|x\|\to +\infty}\inf f_{i_0}(x).
\end{equation}
Because the set $\lbrack f(K) \rbrack^I_{f(x_0)}$ is bounded and $i_0\in I$, we can easily prove that the number $b$ is finite. For each $r>0$, we consider the following quantity:
$$m(r):=\inf_{x\in K^{i_0}_{x_0},\|x\|\geq r}f_{i_0}(x).$$
Then, we can prove that the function $m$ is a nondecreasing function with $\lim_{r\to +\infty}m(r)=b$. So for each $k\in \mathbf{N}$, there exists $r_k>k$ such that for any $r\geq r_k$,
$$m(r)\geq b-\frac{1}{6k}.$$
On  the one hand, it is clear that $m(r)\leq\inf_{x\in K^{i_0}_{x_0},\|x\|\geq r,f(x)\leq f(x_0)}f_{i_0}(x)$ for any $r>0$, since $\{x\in K^{i_0}_{x_0}|\|x\|\geq r,f(x)\leq f(x_0)\}\subseteq\{x\in K^{i_0}_{x_0}|\|x\|\geq r\}$. On the other hand, for any $r>0$ and $y\in \{x\in K^{i_0}_{x_0}|\|x\|\geq r\}$, if $f_{i_0}(y)>f_{i_0}(x_0)$, then we have $f_{i_0}(y)>\inf_{x\in K^{i_0}_{x_0},\|x\|\geq r,f(x)\leq f(x_0)}f_{i_0}(x)$. If $f_{i_0}(y)\leq f_{i_0}(x_0)$, then we have $f_{i_0}(y)\geq\inf_{x\in K^{i_0}_{x_0},\|x\|\geq r,f(x)\leq f(x_0)}f_{i_0}(x)$. Thus, we can deduce that for any $y\in \{x\in K^{i_0}_{x_0}|\|x\|\geq r\}$, $f_{i_0}(y)\geq\inf_{x\in K^{i_0}_{x_0},\|x\|\geq r,f(x)\leq f(x_0)}f_{i_0}(x)$. So, $m(r)\geq\inf_{x\in K^{i_0}_{x_0},\|x\|\geq r,f(x)\leq f(x_0)}f_{i_0}(x)$ for any $r>0$. Therefore,
$$m(r)=\inf_{x\in K^{i_0}_{x_0},\|x\|\geq r,f(x)\leq f(x_0)}f_{i_0}(x), \quad \forall r>0.$$
Hence, we can choose $x_k\in K^{i_0}_{x_0}$ satisfying with $\|x_k\|>2r_k$ such that $f(x_k)\leq f(x_0)$ and $f_{i_0}(x_k)<m(2r_k)+\frac{1}{6k}$. As a consequence,
$$m(2r_k)\leq f_{i_0}(x_k)<m(2r_k)+\frac{1}{6k}\leq b+\frac{1}{6k}\leq m(r_k)+\frac{1}{3k}=\inf_{x\in K^{i_0}_{x_0},\|x\|\geq r_k}f_{i_0}(x)+\frac{1}{3k}.$$
Applying the Ekeland variational principle \cite{Ekland} (see, e.g., \cite[Theorem 2.26]{MK1}, \cite{CLark05}, \cite[Corollary 2.1]{Loridan} and \cite[Theorem 2.1]{Tammer}) to the function $f_{i_0}$ on the closed set $\{x\in K^{i_0}_{x_0}|\|x\|\geq r_k\}$ with the parameter $\epsilon=\frac{1}{3k}$ and $\lambda=\frac{\|x_k\|}{2}$,  we can find $y_k\in K^{i_0}_{x_0}$ with $\|y_k\|\geq r_k$ satisfying the following conditions:
\begin{itemize}
	\item [(a)] $m(r_k)\leq f_{i_0}(y_k)\leq f_{i_0}(x_k)$.
	\item [(b)] $\|y_k-x_k\|\leq \lambda$.
	\item [(c)] $f_{i_0}(y_k)\leq f_{i_0}(x)+\frac{\epsilon}{\lambda}\|x-y_k\|$ for all $x\in K^{i_0}_{x_0}$ with $\|x\|\geq r_k$.
\end{itemize}
It follows from the above condition (a) and $f(x_k)\leq f(x_0)$ that the set $\{f(y_k)\}\subseteq\lbrack f(K) \rbrack_{f(x_0)}$ and $f_{i_0}(y_k)\to b$ as $k\to +\infty$. By the above condition (b), we can deduce
$$r_k<\frac{\|x_k\|}{2}\leq \|y_k\|\leq \frac{3}{2}\|x_k\|.$$
It yields that $\|y_k\|\to +\infty$ as $k\to +\infty$. By the above condition (c), we can know that for each $k\in \mathbf{N}$, $y_k$ is a minimizer of the following  problem:
$$\min_{x\in K^{i_0}_{x_0},\|x\|\geq r_k}f_{i_0}(x)+\frac{\epsilon}{\lambda}\|x-y_k\|.$$
Thus, by $r_k<\|y_k\|$ for all $k$ and the necessary optimality condition given by Lemma \ref{optimality condition} applied to the previous  problem, we can find $\alpha=(\alpha_1,\alpha_2,\dots, \alpha_s)\in \mathbf{R}^s_+$ such that
\begin{equation}\label{08053}
	0\in \alpha_{i_0}\partial\lbrack f_{i_0}(\cdot)+\frac{\epsilon}{\lambda}\|\cdot-y_k\|\rbrack(y_k)+\sum_{i\neq i_0}\alpha_i \partial f_i(y_k)+N(y_k;K),
\end{equation}
$$\alpha_i(f_i(y_k)-f_i(x_0))=0, i=1,2,\dots,i_0-1,i_0+1,\dots,s\quad \mbox{and}\quad \sum_{i=1}^{s}\alpha_i=1.$$
By Lemma \ref{classical subdifferential} and Lemma \ref{sum fule}, we  obtain
$$\alpha_{i_0}\partial\lbrack f_{i_0}(\cdot)+\frac{\epsilon}{\lambda}\|\cdot-y_k\|\rbrack(y_k)\subseteq \alpha_{i_0} \partial f_{i_0}(\cdot)+\alpha_{i_0} \frac{\epsilon}{\lambda}\mathbf{B}^n.$$
The previous inclusion   together with  (\ref{08053}) yields
$$0\in \sum_{i=1}^{s}\alpha_i \partial f_i(y_k)+N(y_k;K)+\alpha_{i_0} \frac{\epsilon}{\lambda}\mathbf{B}^n.$$
From the definition of the extended Rabier function $\nu$ we have
$$\nu(y_k)\leq \alpha_{i_0}\frac{\epsilon}{\lambda}\leq \frac{\epsilon}{\lambda}=\frac{2}{3k\|x_k\|}\leq \frac{1}{\|y_k\|}.$$
It follows that $\nu(y_k)\|y_k\|\leq \frac{1}{k}$ for all $k\in \mathbf{N}$. Therefore, $\nu(y_k)\|y_k\|\to 0$ as $k\to +\infty$. Because the set $\lbrack f(K) \rbrack^I_{f(x_0)}$ is bounded and $\{f(y_k)\}\subseteq\lbrack f(K) \rbrack_{f(x_0)}$, we have that $\{f_I(y_k)\}$ has an accumulation point in $\mathbf{R}^{|I|}$, say $y$. Thus, $y\in K^{I}_{\infty,\leq f(x_0)}(f, K)$, which is a contradiction with assumption $K^{I}_{\infty,\leq f(x_0)}(f, K)=\emptyset$, and so ${\rm(iii)}\Rightarrow {\rm(i)}$.
\par Next, we prove ${\rm(iv)}\Rightarrow {\rm(i)}$. Supposed on the contrary that $f|_{K}$ is not proper with respect to $I$ at the sublevel $f(x_0)$. Then there exists a sequence $\{z_k\}\subseteq K$ such that $f(z_k)\leq f(x_0)$ and $\|z_k\|\to +\infty$ as $k\to +\infty$. Similarly to the above proof, we consider the set $K^{i_0}_{x_0}$  and the number $b$  as defined in  (\ref{0802Kx}) and (\ref{0802b}), respectively. For each $k\in \mathbf{N}$, define the following scalar optimization problem:
$$\min_{x\in K^{i_0}_{x_0}, \|x\|^2=\|z_k\|^2}f_{i_0}(x).$$
Then the constraint set of the above optimization problem is  nonempty and compact, which implies that this problem admits an optimal solution, denoted by $u_k$ for each $k\in \mathbf{N}$. By applying Lemma \ref{optimality condition} to this problem, we can find $(\beta, \eta)=((\beta_1,\beta_2,\dots, \beta_s),\eta)\in \mathbf{R}^s_+\times \mathbf{R}$ such that
$$0\in \sum_{i=1}^{s}\beta_i \partial f_i(u_k)+\eta u_k+N(u_k;K),$$
$$\beta_i(f_i(u_k)-f_i(x_0))=0, i=1,2,\dots,i_0-1,i_0+1,\dots,s\quad \mbox{and}\quad \sum_{i=1}^{s}\beta_i+|\eta|=1.$$
Then, we have $u_k\in \Gamma(f, K)$. Since $\{z_k\}\subseteq K_{x_0}\subseteq K^{i_0}_{x_0}$, we have $z_k\in\{x\in K^{i_0}_{x_0}|\|x\|^2=\|z_k\|^2\}$ for all $k\in \mathbf{N}$. And so $f_{i_0}(u_k)\leq f_{i_0}(z_k)$ for all $k\in \mathbf{N}$. Thus, we  obtain the following properties of the sequence $\{u_k\}$:
\begin{itemize}
	\item [(a')] $\{u_k\}\subseteq \Gamma(f, K)$.
	\item [(b')] $\|u_k\|^2=\|z_k\|^2\to +\infty$ as $k\to +\infty$.
	\item [(c')] $f_{i_0}(u_k)\leq f_{i_0}(z_k)$ for all $k\in \mathbf{N}$.
	\item [(d')] $f(u_k)\leq f(x_0)$ for all $k\in \mathbf{N}$.
\end{itemize}
Because the set $\lbrack f(K) \rbrack^I_{f(x_0)}$ is bounded and $\{f(u_k)\}\subseteq \lbrack f(K) \rbrack_{f(x_0)}$, we have that $\{f_I(u_k)\}$ has an accumulation point in $\mathbf{R}^{|I|}$, say $y'$. Thus, $y'\in T^{I}_{\infty,\leq f(x_0)}(f, K)$, which is a contradiction with assumption $T^{I}_{\infty,\leq f(x_0)}(f, K)=\emptyset$. This yields ${\rm(iv)}\Rightarrow {\rm(i)}$.

Finally, we verify the last statement. Assume that the conclusion ${\rm(i)}$ holds. Since $f$ is  weakly section-bounded from below on $K$ with respect to $I$ and $x_0\in K$, there exist a nonempty index set $I$ and $x_0\in K$ such that the set $\lbrack f(K) \rbrack^I_{f(x_0)}$ is bounded. If the set $\{x\in K|f(x)\leq f(x_0)\}$  is unbounded, then there exists a sequence $\{v_k\}\subseteq \{x\in K|f(x)\leq f(x_0)\}$ such that $\|v_k\|\to +\infty$ as $k\to +\infty$. It follows from conclusion (i) that there exists $i_1\in I$ such that $f_{i_1}(v_k)\to \infty$ as $k\to +\infty$, which is a contradiction with the boundedness of the set $\lbrack f(K) \rbrack^I_{f(x_0)}$. Thus, the set $\{x\in K|f(x)\leq f(x_0)\}$  is bounded and, from the continuity of $f$, we have that the section $\lbrack f(K) \rbrack_{f(x_0)}$ is compact. This completes the proof.\cvd

\begin{remark}\label{Rem:4.5}
	 By the proof of Theorem \ref{relationship}, it follows that  $I$ is an index set such that the set $\lbrack f(K) \rbrack^I_{f(x_0)}$ is bounded. So, if $I=\{1,2,\dots,s\}$ in Theorem \ref{relationship}, then weak section-boundedness from below implies the section-boundedness from below (see also  Remark \ref{0805weak}).
\end{remark}

 Combining Theorem \ref{relationship} with Remark \ref{Rem:4.5}, we recover the following result stated in \cite{Kim3}.

\begin{corollary}\cite[Theorem 3.1]{Kim3}
	 Assume that $f$ is section-bounded from below on $K$. Then the following assertions are equivalent:
    \begin{itemize}
	\item[\rm(i)] $f|_{K}$ is proper at the sublevel $f(x_0)$.
	\item[\rm(ii)] $f|_{K}$ satisfies the Palais-Smale condition at the sublevel $f(x_0)$.
	\item[\rm(iii)] $f|_{K}$ satisfies the weak Palais-Smale condition at the sublevel $f(x_0)$.
    \item[\rm(iv)] $f|_{K}$ satisfies M-tame condition at the sublevel $f(x_0)$.
    \end{itemize}
 Moreover, the set $\{x\in K|f(x)\leq f(x_0)\}$ and the section $\lbrack f(K) \rbrack_{f(x_0)}$ are compact if any of the  conditions (i)-(iv) is fulfilled.
\end{corollary}

\section{Existence of weak Pareto efficient solutions for $\text{VOP}(K,f)$}
\setcounter{equation}{0}

In this section, we  first present a necessary and sufficient condition for the existence of weak Pareto efficient solutions of $\text{VOP}(K,f)$, from which we obtain further particular existence conditions that generalize those obtained in \cite{Kim3}.

Let $y_0\in (\mathbf{R}\cup\{\infty\})^s$ and the nonempty index set $I\subseteq\{1,2,\dots,s\}$, denote the set
$$K^I_{0, \leq y_0}(f, K):=\{f_I(x)\in \mathbf{R}^{|I|}\vert x\in K, f(x)\leq y_0, \nu(x)=0\},$$
where the function $\nu$ is the extended Rabier function. The previous definition, which generalizes the analogous one given in \cite{Kim3}, is motivated by the fact that if ${\rm SOL}^{W}(K,f)$  is nonempty, then there exist a nonempty index set $I\subseteq\{1,2,\dots,s\}$ and $x_0\in {\rm SOL}^{W}(K,f)$ such that $\nu(x_0)=0$, $y_0=f(x_0)$ with $f_I(x)\equiv f_I(x_0)$ for any $x\in K_{x_0}$ and so $f_I(x_0)\in K^I_{0, \leq y_0}(f, K)$ (see the proof of Theorem \ref{gpro}). Notice that if $I=\{1,2,\dots,s\}$, then $K^I_{0, \leq f(x_0)}(f, K)=K_{0, \leq f(x_0)}(f, K)$, where $K_{0, \leq y_0}(f, K):=\{f(x)\in \mathbf{R}^{s}\vert x\in K, f(x)\leq y_0, \nu(x)=0\}$ (see also e.g. \cite{Kim3}).



\begin{theorem}\label{gpro}
    The following assertions are equivalent:
    \begin{itemize}
	\item[\rm(i)] $\rm{VOP}$$(K,f)$ admits a weak Pareto efficient solution.
	\item[\rm(ii)] There exist a nonempty index set $I\subseteq\{1,2,\dots,s\}$ and $x_0\in K$ such that $f$ is  weakly section-bounded from below on $K$ with respect to $I$ and $x_0$, and the inclusion $\widetilde{K}^I_{\infty, \leq f(x_0)}(f, K)\subseteq K^I_{0, \leq f(x_0)}(f, K)$ holds.
	\item[\rm(iii)] There exist a nonempty index set $I\subseteq\{1,2,\dots,s\}$ and $x_0\in K$ such that $f$ is  weakly section-bounded from below on $K$ with respect to $I$ and $x_0$, and the inclusion $K^I_{\infty, \leq f(x_0)}(f, K)\subseteq K^I_{0, \leq f(x_0)}(f, K)$ holds.
    \item[\rm(iv)] There exist a nonempty index set $I\subseteq\{1,2,\dots,s\}$ and $x_0\in K$ such that $f$ is  weakly section-bounded from below on $K$ with respect to $I$ and $x_0$, and the inclusion $T^I_{\infty, \leq f(x_0)}(f, K)\subseteq K^I_{0, \leq f(x_0)}(f, K)$ holds.
    \end{itemize}
\end{theorem}
{\it Proof}\,\,
First, we show ${\rm(i)}\Rightarrow{\rm(ii)},{\rm(i)}\Rightarrow{\rm(iii)}$ and ${\rm(i)}\Rightarrow{\rm(iv)}$. Since ${\rm SOL}^{W}(K, f)\neq\emptyset$, by Proposition \ref{property of nonempty}, there exist $i_0\in \{1,\dots,s\}$ and $x_0\in K$ such that $f_{i_0}(x)\equiv f_{i_0}(x_0)$ for any $x\in K_{x_0}$, where $K_{x_0}:=\{x\in K| f_i(x)\leq f_i(x_0), i=1,2,\dots,s\}$ and $x_0\in {\rm SOL}^{W}(K, f)$. Consider the following set:
$$I:=\{i\in \{1,2,\dots,s \}| f_{i}(x)\equiv f_{i}(x_0),\forall x\in K_{x_0}\}\neq\emptyset.$$
Clearly, $f$ is  weakly section-bounded from below on $K$ with respect to $I$ and $x_0$, since $f_{i}(x)\equiv f_{i}(x_0),\forall x\in K_{x_0}$. Assume that  $I=\{s_1,s_2,\dots,s_q\}$. Then $|I|=q$. By the definition of   $I$,  the set $[f(K)]^I_{f(x_0)}$ is just a singleton $\{f_I(x_0)=(f_{s_1}(x_0),f_{s_2}(x_0),\dots,f_{s_q}(x_0))\}$. It suffices to prove this singleton is contained in the set $T^I_{\infty, \leq f(x_0)}(f, K)$. Since $x_0\in {\rm SOL}^{W}(K, f)$,  for any $x\in K$, there exists $i_x\in \{1,2,\dots,s\}$ such that $f_{i_x}(x)-f_{i_x}(x_0)\geq 0$. Then,
$$\max_{1\leq i\leq s}\{f_i(x)-f_i(x_0)\}\geq 0,\, \forall x\in K.$$
Since $x_0\in K$, then
$$x_0\in \mathop{\arg\min}_{K}\max_{1\leq i\leq s}\{f_i(x)-f_i(x_0)\},$$
and, by the optimality condition of Lemma \ref{optimality condition},  there exists a nonzero $\alpha\in \mathbf{R}$ such that
\begin{equation}\label{0811caculate}
	0\in \alpha\partial \lbrack\max_{1\leq i\leq s}\{f_i(x)-f_i(x_0)\}\rbrack+N(x_0; K).
\end{equation}
Without loss of generality, we assume $\alpha:=1$. By the  Lemma \ref{maximum rule},  there exist $ \lambda_i\geq 0, i\in I(x_0)$ with $\sum_{i\in I(x_0)}\lambda_i=1$ such that
\begin{equation}\label{0811SC}
	\partial \lbrack\max_{1\leq i\leq s}\{f_i(x)-f_i(x_0)\}\rbrack\subseteq\sum_{i\in I(x_0)}\lambda_i\partial f_i(x_0),
\end{equation}
where $I(x_0):=\{i\in \{1,2,\dots,s\}| \max_{1\leq i\leq s}\{f_i(x_0)-f_i(x_0)\}=f_i(x_0)-f_i(x_0)\}$. Clearly, $I(x_0)=\{1,2,\dots,s\}$. Thus, by relationships (\ref{0811caculate}) and (\ref{0811SC}), we obtain that there exists $\lambda\in \mathbf{R}^s_+$ with $\sum^s_{i=1}\lambda_i=1$ such that
$$0\in \sum^s_{i=1}\lambda_i\partial f_i(x_0)+N(x_0; K).$$
By this and recalling the definition of the extended Rabier function $\nu$, we have $\nu(x_0)=0$. Therefore, $f_I(x_0)\in K^I_{0, \leq f(x_0)}(f, K)$. Thus, ${\rm(i)},{\rm(iii)}$ and ${\rm(iv)}$ follow immediately from these facts.
\par Next, we verify ${\rm(ii)}\Rightarrow {\rm(i)},{\rm(iii)}\Rightarrow {\rm(i)}$ and ${\rm(iv)}\Rightarrow {\rm(i)}$. By assumptions, there exist the nonempty index set $I=\{s_1,s_2,\dots,s_q\}\subseteq\{1,2,\dots,s\}$ with $|I|=q$ and $x_0\in K$ such that the set
\begin{equation*}
	\begin{split}
		Y:&=[f(K)]^I_{f(x_0)}=\{f_I\in \mathbf{R}^{|I|}|(f_1,f_2,\dots,f_s)\in \lbrack f(K) \rbrack_{f(x_0)}\}\\
		&=\{f_I(x)\in \mathbf{R}^{|I|}|x\in K, f_i(x)\leq f_i(x_0), i=1,2,\dots,s\}\\
		&=\{(f_{s_1}(x),f_{s_2}(x),\dots,f_{s_q}(x))\in \mathbf{R}^{|I|}|x\in K, f_i(x)\leq f_i(x_0), i=1,2,\dots,s\}\subseteq\mathbf{R}^q
	\end{split}
\end{equation*}
is bounded, so that the closure $\bar Y$ of the set $Y$ is a nonempty, compact set in $\mathbf{R}^q$. By  Claim 1 in \cite[Theorem 4.1]{Kim3}, we have that there exists $\alpha\in \mbox{ int } \mathbf{R}^{q}_+$ such that the scalar optimization problem:
\begin{equation}\label{0802op}
	\min_{y\in \bar Y}\langle \alpha, y \rangle
\end{equation}
has the unique optimal solution, denoted by $\bar y:=(\bar y_{s_1},\bar y_{s_2},\dots,\bar y_{s_q})\in \bar Y$. Set $K_{x_0}:=\{x\in K| f_i(x)\leq f_i(x_0), i=1,2,\dots,s\}$. Then $K_{x_0}$ is nonempty. The proofs consists of three steps:

Step 1. We claim that if $K_{x_0}\bigcap f^{-1}_I(\bar y)\neq\emptyset$, then for each $\bar x\in K_{x_0}\bigcap f^{-1}_I(\bar y)$, $\bar x$ is a weak Pareto efficient solution of $\rm{VOP}$$(K,f)$.

Suppose on the contrary that there exists $x'\in K$ such that $f(x')<f(\bar x)$.  Since $\bar x\in K_{x_0}$ then $x'\in K_{x_0}$, so that  $f_I(x')\in Y\subseteq \bar Y$. From $\bar x\in K_{x_0}\bigcap f^{-1}_I(\bar y)$ it follows that $\langle \alpha, f_I(x') \rangle<\langle \alpha, f_I(\bar x) \rangle=\min_{y\in \bar Y}\langle \alpha, y \rangle$, which is a contradiction.
\par Step 2. We claim that if $K_{x_0}\bigcap f^{-1}_I(\bar y)=\emptyset$, then $\bar y\in K^I_{\infty, \leq f(x_0)}(f, K)\bigcap T^I_{\infty, \leq f(x_0)}(f, K)\bigcap\widetilde{K}^{I}_{\infty,\leq f(x_0)}(f, K)$.

Indeed, let $f(x_0)=y^0=(y_1^0, y_2^0, \dots,y_s^0)$ and $f_I(x_0)=y^0_I=(y_{s_1}^0, y_{s_2}^0, \dots,y_{s_q}^0)$. Since $\bar y\in \bar Y$, we have $\bar y_i\leq y^0_i$ for all $i\in I$. Consider the index set
$$J:=\{i\in I|\bar y_i<y^0_i\}.$$
Then $J\subseteq I$ is nonempty:  in fact, if, on the contrary, $\bar y=y^0_I$, then $x_0\in K_{x_0}\bigcap f^{-1}_I(\bar y)\neq\emptyset$, which  contradicts $K_{x_0}\bigcap f^{-1}_I(\bar y)=\emptyset$.
We first show that
\begin{equation}\label{08021}
	\arg\min\{\sum_{i\in J\subseteq I}\alpha_i y_i|y\in \bar Y\}=\{\bar y\}.
\end{equation}
In fact, by compactness of $\bar Y$, $\arg\min\{\sum_{i\in J\subseteq I}\lambda_i y_i|y\in \bar Y\}\neq\emptyset$.
Let $y\in \arg\min\{\sum_{i\in J\subseteq I}\alpha_i y_i|y\in \bar Y\}$. Then we have
$$\langle \alpha, y\rangle=\sum_{i\in J}\alpha_i y_i+\sum_{i\in I\backslash J}\alpha_i y_i\leq \sum_{i\in J}\alpha_i \bar y_i+\sum_{i\in I\backslash J}\alpha_i y^0_i=\sum_{i\in J}\alpha_i \bar y_i+\sum_{i\in I\backslash J}\alpha_i \bar y_i=\langle \alpha, \bar y\rangle.$$
Since $\bar y$ is the unique minimizer of the scalar optimization problem (\ref{0802op}), then $y=\bar y$.

Define the function $\varphi: \mathbf{R}^n\to \mathbf{R}$ by
$$\varphi(x):=\sum_{i\in J\subseteq I}\alpha_i f_i(x).$$
We next show that
\begin{equation}\label{08022}
	\min\{\sum_{i\in J\subseteq I}\alpha_i y_i|y\in \bar Y\}=\inf_{x\in K_{x_0}}\varphi(x).
\end{equation}
Indeed, for any $x\in K_{x_0}$, we have $f(x)\in f(K_{x_0})$. So $f_I(x)\in Y$ and so
\begin{equation}\label{08041}
	\min\{\sum_{i\in J\subseteq I}\alpha_i y_i|y\in \bar Y\}\leq \sum_{i\in J\subseteq I}\alpha_i f_i(x)=\varphi(x).
\end{equation}
This implies that
$$\min\{\sum_{i\in J\subseteq I}\alpha_i y_i|y\in \bar Y\}\leq \inf_{x\in K_{x_0}}\varphi(x).$$
Since $\bar y\in \bar Y$, there exists a sequence $\{x_k\}\subseteq K_{x_0}$ such that $f_I(x_k)\to \bar y$ as $k\to +\infty$. Then for each $k\in \mathbf{N}$, we have
$$\inf_{x\in K_{x_0}}\varphi(x)\leq \varphi(x_k)=\sum_{i\in J\subseteq I}\alpha_i f_i(x_k).$$
By  (\ref{08021}) and letting $k\to +\infty$, we obtain
\begin{equation}\label{08042}
	\inf_{x\in K_{x_0}}\varphi(x)\leq \sum_{i\in J\subseteq I}\alpha_i \bar y_i=\min\{\sum_{i\in J\subseteq I}\alpha_i y_i|y\in \bar Y\}.
\end{equation}
Combining the inequalities  (\ref{08041}) and  (\ref{08042}), we  obtain equality (\ref{08022}).

If $\varphi(x)$ attains its infimum on $K_{x_0}$, say $ x^*\in K_{x_0}$, then by the equalities (\ref{08021}) and (\ref{08022}), we have $\bar y=f_I(x^*)\in Y\subseteq\bar Y$ and so $x^*\in K_{x_0}\bigcap f^{-1}_I(\bar y)\neq\emptyset$, which is a contradiction. Therefore, $\varphi(x)$ does not attains its infimum on $K_{x_0}$ and since $K_{x_0}$ is closed, we have that $K_{x_0}$ is unbounded. Moreover,  let $\epsilon=\frac{1}{3k},k\in \mathbf{N}$, there exists a sequence $\{x_k\}\subseteq K_{x_0}$ with $\|x_k\|\to +\infty $ as $k\to +\infty$, such that
$$\inf_{x\in K_{x_0}}\varphi(x)\leq \varphi(x_k)\leq \inf_{x\in K_{x_0}}\varphi(x)+\epsilon.$$
Then $\varphi(x_k)\to \inf_{x\in K_{x_0}}\varphi(x)$ as $k\to +\infty$. Let $\lambda=\frac{\|x_k\|}{2}$. By using the Ekeland variational principle as in the proof of Theorem \ref{relationship} applied  to the function $\varphi$ on the closed set $K_{x_0}$ with  parameters $\epsilon=\frac{1}{3k}$ and $\lambda=\frac{\|x_k\|}{2}$,  we can find a sequence $\{u_k\}\subseteq K_{x_0}$ satisfying the following conditions:
\begin{itemize}
	\item [(a$^1$)] $\inf_{x\in K_{x_0}}\varphi(x)\leq \varphi(u_k)\leq \varphi(x_k)$.
	\item [(b$^1$)] $\|u_k-x_k\|\leq \lambda$.
	\item [(c$^1$)] $\varphi(u_k)\leq \varphi(x)+\frac{\epsilon}{\lambda}\|x-u_k\|$ for all $x\in K_{x_0}$.
\end{itemize}
Then, by the above condition (a$^1$), we have $\varphi(u_k)\to \inf_{x\in K_{x_0}}\varphi(x)$. By the above condition (b$^1$), we can deduce
$$\frac{\|x_k\|}{2}\leq \|u_k\|\leq \frac{3}{2}\|x_k\|.$$
It yields $\|u_k\|\to +\infty$ as $k\to +\infty$. By the above condition (c$^1$), we have that, for each $k\in \mathbf{N}$, $u_k$ is a minimizer of the following scalar optimization problem:
$$\min_{x\in K_{x_0}}\varphi(x)+\frac{\epsilon}{\lambda}\|x-u_k\|.$$
Thus, by applying Lemma \ref{optimality condition} to this problem,  for each $k$, there exists a pair $(\kappa, \beta)\in \mathbf{R}_+\times\mathbf{R}^s_+$ satisfying the following relationships:
\begin{equation}\label{08043}
	0\in \kappa\partial\lbrack \varphi(\cdot)+\frac{\epsilon}{\lambda}\|\cdot-u_k\|\rbrack(u_k)+\sum_{i=1}^s\beta_i \partial f_i(u_k)+N(u_k;K),
\end{equation}
$$\beta_i(f_i(u_k)-f_i(x_0))=0, i=1,2,\dots,s\quad \mbox{and}\quad \kappa+\sum_{i=1}^{s}\beta_i=1.$$
Note that from Lemma \ref{sum fule} and Lemma \ref{classical subdifferential}, we  obtain
\begin{equation}\label{08044}
	\partial\lbrack \varphi(\cdot)+\frac{\epsilon}{\lambda}\|\cdot-u_k\|\rbrack(u_k)\subseteq \partial \varphi(u_k)+\frac{\epsilon}{\lambda}\mathbf{B}^n\subseteq \sum_{i\in J\subseteq I}\alpha_i \partial f_i(u_k)+\frac{\epsilon}{\lambda}\mathbf{B}^n.
\end{equation}
It follows from the boundedness of the set $[f(K)]^I_{f(x_0)}$ and $\{f(u_k)\}\subseteq \lbrack f(K) \rbrack_{f(x_0)}$, we have that $\{f_I(u_k)\}$ has an accumulation point in $\mathbf{R}^{|I|}$, say $y$. It follows from $\varphi(u_k)\to \inf_{x\in K_{x_0}}\varphi(x)$ as $k\to +\infty$,  (\ref{08021}) and (\ref{08022}) that we can get $y=\bar y$. Thus, we have $f_I(u_k)\to \bar y$ as $k\to +\infty$. By the definition of the index set $J\subseteq I$, for all $i\in J$ and all $k$ sufficiently large, we have $f_i(u_k)-y^0_i<0$, that is, $f_i(u_k)-f_i(x_0)<0$. Thus, $\beta_i=0$ for all $i\in J$ and all $k$ sufficiently large. By the relationships (\ref{08043}) and (\ref{08044}), we have
$$0\in \sum_{i\in J\subseteq I}\kappa\alpha_i \partial f_i(u_k)+\sum_{i\notin J}\beta_i\partial f_i(u_k)+N(u_k;K)+\kappa\frac{\epsilon}{\lambda}\mathbf{B}^n,$$
for all $k$ sufficiently large. By the definition of the extended Rabier function $\nu$ and $0<\frac{\epsilon}{\lambda}\leq \frac{1}{k\|u_k\|}$, we  obtain that
$$\nu(u_k)\leq \frac{\kappa}{\sum_{i\in J}\kappa\alpha_i+\sum_{i\notin J}\beta_i}\frac{\epsilon}{\lambda}\leq \frac{1}{\sum_{i\in J}\alpha_i}\frac{\epsilon}{\lambda}\leq \frac{1}{\sum_{i\in J}\alpha_i}\frac{1}{k\|u_k\|}.$$
Thus, we get the estimate
 $$\|u_k\|\nu(u_k)\leq \frac{1}{\sum_{i\in J}\alpha_i}\frac{1}{k},$$
and so $\|u_k\|\nu(u_k)\to 0$ as $k\to +\infty$. Therefore, $\bar y\in K^{I}_{\infty,\leq f_{i_0}(x_0)}(f, K)$.
\par Next, we prove $\bar y\in T^{I}_{\infty,\leq f(x_0)}(f, K)$. For each $k\in \mathbf{N}$, we consider the following scalar optimization problem:
$$\min_{x\in K_{x_0}, \|x\|^2=\|x_k\|^2}\varphi(x),$$
where $\{x_k\}\subseteq K_{x_0}$, $\varphi(x_k)\to \inf_{x\in K_{x_0}}\varphi(x)$ and $\|x_k\|\to +\infty$ as $k\to +\infty$. Since the constraint set is nonempty and compact,  this problem admits an optimal solution denoted by $v_k$ and, by  Lemma \ref{optimality condition}, we can find $(\beta, \eta)=(\beta_0, \beta_1,\beta_2,\dots, \beta_s, \eta)\in \mathbf{R}^{s+1}_+\times \mathbf{R}$ such that
\begin{equation}\label{08045}
	0\in \beta_0\partial \varphi(v_k)+\sum_{i=1}^{s}\beta_i\partial f_i(v_k)+\eta v_k+N(v_k;K),
\end{equation}
$$\beta_i(f_i(v_k)-f_i(x_0))=0, i=1,2,\dots,s\quad \mbox{and}\quad \sum_{i=0}^{s}\beta_i+|\eta|=1.$$
Note that the subdifferential sum rule from Lemma \ref{sum fule} gives us
\begin{equation}\label{08046}
	\partial \varphi(v_k)\subseteq \sum_{i\in J\subseteq I}\alpha_i \partial f_i(v_k).
\end{equation}
On the other hand,  from the boundedness of the set $[f(K)]^I_{f(x_0)}$ and $\{f(v_k)\}\subseteq \lbrack f(K) \rbrack_{f(x_0)}$, we have that $\{f_I(v_k)\}$ has an accumulation point in $\mathbf{R}^{|I|}$, say $y'$. Then, by $\inf_{x\in K_{x_0}}\varphi(x)\leq \varphi(v_k)\leq \varphi(x_k)$ for all $k\in \mathbf{N}$ and $\varphi(x_k)\to \inf_{x\in K_{x_0}}\varphi(x)$ and $\|v_k\|=\|x_k\|\to +\infty$ as $k\to +\infty$, we have $\varphi(v_k)\to \inf_{x\in K_{x_0}}\varphi(x)$ as $k\to +\infty$, and combining with the equations (\ref{08021}) and (\ref{08022}), we can get $y'=\bar y$. Thus, we have $f_I(v_k)\to \bar y$ as $k\to +\infty$. By the definition of the index set $J\subseteq I$, for all $i\in J$ and all $k$ sufficiently large, we have $f_i(v_k)-y^0_i<0$, that is, $f_i(u_k)-f_i(x_0)<0$. So $\beta_i=0$ for all $i\in J$ and all $k$ sufficiently large. Thus, by relationships (\ref{08045}) and (\ref{08046}), we have
$$0\in \sum_{i\in J\subseteq I}\beta_0\alpha_i \partial f_i(v_k)+\sum_{i\notin J}\beta_i\partial f_i(v_k)+\eta v_k+N(v_k;K),$$
for all $k$ sufficiently large. Then, for all $k$ sufficiently large, we can easily prove $v_k\in \Gamma(f, K)$. So, according to the above discussion, we obtain the following properties of the sequence $\{v_k\}$:
\begin{itemize}
	\item [(a$^2$)] $\{v_k\}\subseteq \Gamma(f, K)\bigcap K_{x_0}$.
	\item [(b$^2$)] $\|v_k\|\to +\infty$ as $k\to +\infty$.
	\item [(c$^2$)] $f_I(v_k)\to \bar y$ as $k\to +\infty$.
\end{itemize}
By definition, this gives us that $\bar y\in T^{I}_{\infty,\leq f(x_0)}(f, K)$. Thus, we have proved that
$$\bar y\in K^I_{\infty, \leq f(x_0)}(f, K)\bigcap T^I_{\infty, \leq f(x_0)}(f, K).$$
Moreover, by $K^{I}_{\infty,\leq y_0}(f, K)\subseteq \widetilde{K}^{I}_{\infty,\leq y_0}(f, K)$, we have $$\bar y\in K^I_{\infty, \leq f(x_0)}(f, K)\bigcap T^I_{\infty, \leq f(x_0)}(f, K)\bigcap\widetilde{K}^{I}_{\infty,\leq f(x_0)}(f, K).$$

Step 3. Now, we will finish the proof of Theorem \ref{gpro}. In fact, if $K_{x_0}\bigcap f^{-1}_I(\bar y)=\emptyset$, then $$\bar y\in K^I_{\infty, \leq f(x_0)}(f, K)\bigcap T^I_{\infty, \leq f(x_0)}(f, K)\bigcap\widetilde{K}^{I}_{\infty,\leq f(x_0)}(f, K),$$
which yields that $\bar y\in K^I_{0, \leq f(x_0)}(f, K)$,  taking into account the assumptions. Hence, for some $x'\in K_{x_0}$, we have $\bar y=f_I(x')$, which is a contradiction with $K_{x_0}\bigcap f^{-1}_I(\bar y)=\emptyset$. Thus, $K_{x_0}\bigcap f^{-1}_I(\bar y)\neq\emptyset$. By Step 1, we have ${\rm SOL}^{W}(K, f)\neq\emptyset$. This completes the proof.\cvd

\begin{remark}
	In particular, if $I=\{1,2,\dots,s\}$ in Theorem \ref{gpro}, then $K^I_{0, \leq f(x_0)}(f, K)=K_{0, \leq f(x_0)}(f, K)$, $\widetilde{K}^I_{\infty, \leq f(x_0)}(f, K)=\widetilde{K}_{\infty, \leq f(x_0)}(f, K)$, $K^I_{\infty, \leq f(x_0)}(f, K)=K_{\infty, \leq f(x_0)}(f, K)$, $T^I_{\infty, \leq f(x_0)}(f, K)=T_{\infty, \leq f(x_0)}(f, K)$ and weakly section-bounded from below implies the section-bounded from below. Thus, Theorem \ref{gpro} reduces to  \cite[Theorem 4.1]{Kim3}.
\end{remark}

As a direct consequence of Theorem \ref{relationship} and Theorem \ref{gpro}, we  obtain the following result.

\begin{corollary}\label{kong}
	Assume that there exist a nonempty index set $I\subseteq\{1,2,\dots,s\}$ and $x_0\in K$ such that $f$ is  weakly section-bounded from below on $K$ with respect to $I$ and $x_0$. Then $\rm{VOP}$$(K,f)$ admits a weak Pareto efficient solution provided that one of the following equivalent conditions holds:
	\begin{itemize}
	\item[\rm(i)] $f|_{K}$ is proper with respect to $I$ at the sublevel $f(x_0)$.
	\item[\rm(ii)] $f|_{K}$ satisfies the Palais-Smale condition with respect to $I$ at the sublevel $f(x_0)$.
	\item[\rm(iii)] $f|_{K}$ satisfies the weak Palais-Smale condition with respect to $I$ at the sublevel $f(x_0)$.
    \item[\rm(iv)] $f|_{K}$ satisfies M-tame condition with respect to $I$ at the sublevel $f(x_0)$.
    \end{itemize}
\end{corollary}

\begin{remark}
	In particular, if $I=\{1,2,\dots,s\}$ in Corollary  \ref{kong}, then Corollary  \ref{kong} reduces to  \cite[Corollary 4.1]{Kim3}, that guarantees the existence of a Pareto efficient solution for \rm{VOP}$(K,f)$.
\end{remark}

 Corollary \ref{C1} shows that there is a close relationship between the nonemptiness of  weak Pareto efficient solution for VOP$(K, f)$ and nonemptiness of Pareto efficient solution for VOP$(K_{x_0}, f_I)$,
for some nonempty index set $I\subseteq\{1,2,\dots,s\}$ and $x_0\in K$. Then, by \cite[Theorem 4.1]{Kim3} and Corollary \ref{C1}, we can prove the following existence theorem.

\begin{theorem}\label{20231227ex}
    The following assertions are equivalent:
    \begin{itemize}
	\item[\rm(i)] $\rm{VOP}$$(K,f)$ admits a weak Pareto efficient solution.
	\item[\rm(ii)] There exist a nonempty index set $I\subseteq\{1,2,\dots,s\}$ and $x_0\in K$ such that $f$ is  weakly section-bounded from below on $K$ with respect to $I$ and $x_0$, and the inclusion $\widetilde{K}_{\infty, \leq f_I(x_0)}(f_I, K_{x_0})\subseteq K_{0, \leq f_I(x_0)}(f_I, K_{x_0})$ holds.
	\item[\rm(iii)] There exist a  nonempty index set $I\subseteq\{1,2,\dots,s\}$ and $x_0\in K$ such that $f$ is  weakly section-bounded from below on $K$ with respect to $I$ and $x_0$, and the inclusion $K_{\infty, \leq f_I(x_0)}(f_I, K_{x_0})\subseteq K_{0, \leq f_I(x_0)}(f_I, K_{x_0})$ holds.
    \item[\rm(iv)] There exist a nonempty index set $I\subseteq\{1,2,\dots,s\}$ and $x_0\in K$ such that $f$ is  weakly section-bounded from below on $K$ with respect to $I$ and $x_0$, and the inclusion $T_{\infty, \leq f_I(x_0)}(f_I, K_{x_0})\subseteq K_{0, \leq f_I(x_0)}(f_I, K_{x_0})$ holds.
    \end{itemize}
\end{theorem}
{\it Proof}\,\,
As observed in Remark \ref{0805weak}, $f$ is  weakly section-bounded from below on $K$ with respect to $I$ and $x_0$, if and only if $f_I$ is  section-bounded from below on $K_{x_0}$ at $x_0$, that is, the section $[f_I(K_{x_0})]_{f_I(x_0)}$ is bounded.

``$\Rightarrow$'': By Corollary \ref{C1}, if ${\rm SOL}^{W}(K, f)\neq\emptyset$, then
there exists $x_0\in K$ such that the index set $I$ as in \eqref{L1} is nonempty and $x_0\in {\rm SOL}^{S}(K_{x_0}, f_I)$. So, applying \cite[Theorem 4.1]{Kim3} to the problem  $\text{VOP}(K_{x_0},f_I)$, we have that the section $\lbrack f_I(K_{x_0}) \rbrack_{f_I(x_0)}$ is bounded,
 or equivalently,
$f$ is  weakly section-bounded from below on $K$ with respect to $I$ and $x_0$, and  that all the inclusion relationships  of the conditions $(\rm ii)-(\rm iv)$ hold. Hence, we prove $(\rm i)\Rightarrow (\rm ii)$, $(\rm i)\Rightarrow (\rm iii)$ and $(\rm i)\Rightarrow (\rm iv)$ .

``$\Leftarrow$'':Conversely, if there exist a nonempty index set $I\subseteq\{1,2,\dots,s\}$ and $x_0\in K$ such that $f$ is  weakly section-bounded from below on $K$ with respect to $I$ and $x_0$, then
by Remark \ref{0805weak},  we have 
that the section $[f_I(K_{x_0})]_{f_I(x_0)}$ is bounded. If one of the inclusion conditions in $(\rm ii)-(\rm iv)$ holds, then applying \cite[Theorem 4.1]{Kim3} to $\text{VOP}(K_{x_0},f_I)$, we have ${\rm SOL}^{S}(K_{x_0}, f_I)\neq\emptyset$, and so ${\rm SOL}^{W}(K_{x_0}, f_I)\neq\emptyset$ by ${\rm SOL}^{S}(K_{x_0}, f_I)\subseteq{\rm SOL}^{W}(K_{x_0}, f_I)$. By Lemma  \ref{2023801}, we can get ${\rm {SOL}}^{W}(K, f_I)\neq\emptyset$. Let $\bar x\in {\rm SOL}^{W}(K, f_I)$. Then for any $x\in K$, there exists $i\in I\subseteq \{1,2,\dots,s\}$ such that $f_i(x)-f_i(\bar x)\geq 0$.
It follows that $\bar x\in {\rm SOL}^{W}(K, f)$. Hence, ${\rm SOL}^{W}(K, f)\neq\emptyset$ and the proof is completed.
\cvd

\begin{remark}
We notice that  the assumptions of Theorem \ref{20231227ex} are different from those of Theorem \ref{gpro}.  In fact, by the definition, in  Theorem \ref{gpro}  the extended Rabier function and the tangency variety depend on the mapping $f$ and on the set $K$.
 Instead,  in Theorem \ref{20231227ex}, since we replace $K$ with $K_{x_0}$ and $f$ with $f_I$, then the extended Rabier function $\nu$ is changed by

$$\nu_{I,K_{x_0}}(x):=\inf\{\|\sum_{i\in I}\alpha_i v_i+\omega\|\vert v_i\in \partial f_i(x),i\in I, \omega\in N(x;K_{x_0}), \alpha\in \mathbf{R}^{|I|}_+, \sum_{i\in I}\alpha_i=1\},$$
and $\Gamma(f, K)$ is changed by $\Gamma(f_I, K_{x_0})$.

Hence, the following two examples show that we can not guarantee that the equalities $\widetilde{K}^{I}_{\infty,\leq f(x_0)}(f, K)=\widetilde{K}_{\infty,\leq f_I(x_0)}(f_I, K_{x_0})$, $K^{I}_{\infty,\leq f(x_0)}(f, K)= K_{\infty,\leq f_I(x_0)}(f_I, K_{x_0})$, $T^{I}_{\infty,\leq f(x_0)}(f, K)= T_{\infty,\leq f_I(x_0)}(f_I, K_{x_0})$\\ and $K^I_{0, \leq f(x_0)}(f, K)=K_{0, \leq f_I(x_0)}(f_I, K_{x_0})$ hold.
First, consider the following two extended Rabier functions
$$\nu(x):=\inf\{\|\sum_{i=1}^{2}\alpha_i \nabla f_i(x)+\omega\| \vert  \omega\in N(x;K), \alpha\in \mathbf{R}^{2}_+, \sum_{i=1}^{2}\alpha_i=1\}$$
and
$$\nu_{I,K_{x_0}}(x):=\inf\{\|\sum_{i\in I}\alpha_i \nabla f_i(x)+\omega\|\vert \omega\in N(x;K_{x_0}), \alpha\in \mathbf{R}^{|I|}_+, \sum_{i\in I}\alpha_i=1\}$$
 with $I\subseteq\{1,2,\dots,s\}$.
Consider the following two tangency varieties
\begin{equation*}
\begin{split}
	&\Gamma(f, K):=\{x\in K_{x_0}\vert \exists (\alpha, \mu)\in \mathbf{R}^2_+\times \mathbf{R} \mbox{ with }  \sum_{i=1}^{2}\alpha_i+|\mu|=1 \mbox{ such that } \\
		&0\in  \sum_{i=1}^{2}\alpha_i \nabla f_i(x)+\mu x+N(x; K)\}.
\end{split}
		\end{equation*}
and
\begin{equation*}
\begin{split}
	&\Gamma(f_I, K_{x_0}):=\{x\in K_{x_0}\vert  \exists (\alpha, \mu)\in \mathbf{R}_+\times \mathbf{R} \mbox{ with }  \sum_{i\in I}\alpha_i+|\mu|=1 \mbox{ such that } \\
		&0\in  \sum_{i\in I}\alpha_i \nabla f_i(x)+\mu x+N(x; K_{x_0})\}
\end{split}
		\end{equation*}
 with $I\subseteq\{1,2,\dots,s\}$.
\end{remark}


\begin{example}
	Let $f:=(f_1, f_2): \mathbf{R}\to \mathbf{R}^2$ with $f_1(x)=\sin x$, $f_2(x)=\sin x+\cos x$, $K:=\{x\in\mathbf{R}\vert x\in [-\pi+2k\pi, \frac{\pi}{2}+2k\pi], k\in \mathbf{N}\}$ with $I=\{1\}$ and $x_0=\frac{\pi}{4}$.
	

	 First, we  prove that the set $K$ is closed. We rewrite the set $K:=\{x\in\mathbf{R}\vert x\notin (-\frac{3\pi}{2}+2k\pi, -\pi+2k\pi), k\in \mathbf{N}\}$ as $K:=\{x\in\mathbf{R}\vert x\in [-\pi+2k\pi, \frac{\pi}{2}+2k\pi], k\in \mathbf{N}\}=\mathbf{R}\backslash K^*$, where $K^*:=\{x\in\mathbf{R}\vert x\in (-\frac{3\pi}{2}+2k\pi, -\pi+2k\pi), k\in \mathbf{N}\}$. Observe that $K^*=\bigcup_{k\in \mathbf{N}} (-\frac{3\pi}{2}+2k\pi, -\pi+2k\pi)$, and since $(-\frac{3\pi}{2}+2k\pi, -\pi+2k\pi)$ is open for each $k\in \mathbf{N}$, then $K$ is closed.  On the one hand, $f:=(f_{1}, f_{2})$ is  a locally Lipschitz mapping and $f$ is  weakly section-bounded from below on the set $K$ with respect to $I=\{1\}$ and $x_0=\frac{\pi}{4}$. By definitions, it is easy to calculate that the inclusions $\widetilde{K}^I_{\infty, \leq f(x_0)}(f, K)\subseteq \{[-1, -\frac{\sqrt 2}{2}], \frac{\sqrt 2}{2}\}\subseteq K^I_{0, \leq f(x_0)}(f, K)$, $K^I_{\infty, \leq f(x_0)}(f, K)\subseteq \{[-1, -\frac{\sqrt 2}{2}], \frac{\sqrt 2}{2}\}\subseteq K^I_{0, \leq f(x_0)}(f, K)$, $\widetilde{K}_{\infty, \leq f_I(x_0)}(f_I, K_{x_0})\subseteq \{-1\}\subseteq K_{0, \leq f_I(x_0)}(f_I, K_{x_0})$ and $K_{\infty, \leq f_I(x_0)}(f_I, K_{x_0})\subseteq \{-1\}\subseteq K_{0, \leq f_I(x_0)}(f_I, K_{x_0})$ hold. On the other hand, by the definition of the extended Rabier function, let $y_k=-\frac{3\pi}{4}+2k\pi$, $k\in \mathbf{N}$. Then $f_I(y_k)\to -\frac{\sqrt{2}}{2}$ as $k\to +\infty$ and  $\nu(y_k)=0$, $k\in \mathbf{N}$. It follows from $\{y_k\}\subseteq K$,  $\|y_k\|\to +\infty$ as $k\to +\infty$, $f(y_k)\leq f(x_0),\nu(y_k)\equiv 0$ with $k\in \mathbf{N}$ and the definitions of $\widetilde{K}^{I}_{\infty,\leq f(x_0)}(f, K)$, $K^{I}_{\infty,\leq f(x_0)}(f, K)$ that $-\frac{\sqrt{2}}{2}\in \widetilde{K}^{I}_{\infty,\leq f(x_0)}(f, K)$, $-\frac{\sqrt{2}}{2}\in K^{I}_{\infty,\leq f(x_0)}(f, K)$ and $-\frac{\sqrt{2}}{2}\in K^I_{0, \leq f(x_0)}(f, K)$. However, $-\frac{\sqrt{2}}{2}\notin \widetilde{K}_{\infty,\leq f_I(x_0)}(f_I, K_{x_0})$, $-\frac{\sqrt{2}}{2}\notin K_{\infty,\leq f_I(x_0)}(f_I, K_{x_0})$ and $-\frac{\sqrt{2}}{2}\notin K_{0, \leq f_I(x_0)}(f_I, K_{x_0})$. Hence, the three equalities $\widetilde{K}^{I}_{\infty,\leq f(x_0)}(f, K)=\widetilde{K}_{\infty,\leq f_I(x_0)}(f_I, K_{x_0})$, $K^{I}_{\infty,\leq f(x_0)}(f, K)= K_{\infty,\leq f_I(x_0)}(f_I, K_{x_0})$ and $K^I_{0, \leq f(x_0)}(f, K)=K_{0, \leq f_I(x_0)}(f_I, K_{x_0})$  do not hold.
\end{example}

\begin{example}
	Let $f:=(f_1, f_2): \mathbf{R}^2\to \mathbf{R}^2$ with $f_1(x)=x_1x_2$, $f_2(x)=\sin x_1+\cos x_1$, $K:=\{x\in\mathbf{R}^2\vert x_1\geq 0, x_2\geq 0\}$ and $x_0=(x_1, x_2)=(\frac{\pi}{4}, \frac{4}{\pi})$.
 Clearly, the set $K_{x_0}$ is closed. On the one hand, $f:=(f_{1}, f_{2})$ is  a locally Lipschitz mapping and $f$ is  weakly section-bounded from below on $K$ with respect to $I=\{1\}$ and $x_0=(\frac{\pi}{4}, \frac{4}{\pi})$. By definitions, it is easy to calculate that the inclusions $T^I_{\infty, \leq f(x_0)}(f, K)\subseteq [0, 1]= K^I_{0, \leq f(x_0)}(f, K)$ and $T_{\infty, \leq f_I(x_0)}(f_I, K_{x_0})\subseteq \{0, 1\}\subseteq K_{0, \leq f_I(x_0)}(f_I, K_{x_0})$ hold. On the other hand, by the definition of the tangency variety, let $x_k=(\frac{\pi}{4}+2k\pi, \frac{1}{\frac{\pi}{2}+4k\pi})$, $k\in \mathbf{N}$, we can know $f_I(x_k)\to \frac{1}{2}$ as $k\to +\infty$ and $x_k\in \Gamma(f, K)$, $k\in \mathbf{N}$. It follows from $\{x_k\}\subseteq K$,  $f(x_k)\leq f(x_0)$ with $k\in \mathbf{N}$, $\|x_k\|\to +\infty$ as $k\to +\infty$, and the definition of $T^{I}_{\infty,\leq f(x_0)}(f, K)$ that $\frac{1}{2}\in T^{I}_{\infty,\leq f(x_0)}(f, K)$ and $\frac{1}{2}\in K^I_{0, \leq f(x_0)}(f, K)$. However, $\frac{1}{2}\notin T_{\infty,\leq f_I(x_0)}(f_I, K_{x_0})$ and $\frac{1}{2}\notin K_{0, \leq f_I(x_0)}(f_I, K_{x_0})$. Hence, the two equalities $T^{I}_{\infty,\leq f(x_0)}(f, K)= T_{\infty,\leq f_I(x_0)}(f_I, K_{x_0})$ and $K^I_{0, \leq f(x_0)}(f, K)=K_{0, \leq f_I(x_0)}(f_I, K_{x_0})$ are not true.
\end{example}


\begin{remark}
In particular, if $I=\{1,2,\dots,s\}$ in Theorem \ref{20231227ex}, then $[f(K_{x_0})]_{f(x_0)}=[f_I(K_{x_0})]_{f_I(x_0)}$ and
 $[f_I(K_{x_0})]_{f_I(x_0)}$ is weakly bounded with respect to $I$ and $x_0$  if and only if the section $[f(K)]_{f(x_0)}$ is bounded (see Remark \ref{Rem:1}). Applying \cite[Theorem 4.1]{Kim3} to $\rm{VOP} $$(K, f)$, we can prove the following characterization of the nonemptiness of ${\rm{SOL}}^S $$(K, f)$.
\end{remark}

\begin{corollary}\label{20231227co}
The following assertions are equivalent:
    \begin{itemize}
    \item[\rm(i)] $\rm{VOP}$$(K,f)$ admits a Pareto efficient solution.
	\item[\rm(ii)] There exists $x_0\in K$ such that the section $[f(K)]_{f(x_0)}$ is bounded, and the inclusion $\widetilde{K}_{\infty, \leq f(x_0)}(f, K_{x_0})\subseteq K_{0, \leq f(x_0)}(f, K_{x_0})$ holds.
	\item[\rm(iii)] There exists $x_0\in K$ such that the section $[f(K)]_{f(x_0)}$ is bounded, and the inclusion $K_{\infty, \leq f(x_0)}(f, K_{x_0})\subseteq K_{0, \leq f(x_0)}(f, K_{x_0})$ holds.
    \item[\rm(iv)] There exists $x_0\in K$ such that the section $[f(K)]_{f(x_0)}$ is bounded, and the inclusion $T_{\infty, \leq f(x_0)}(f, K_{x_0})\subseteq K_{0, \leq f(x_0)}(f, K_{x_0})$ holds.
    \end{itemize}
\end{corollary}
{\it Proof}\,\,
Let $x_0\in {\rm SOL}^{S}(K, f)$. Then it is easy to prove $x_0\in {\rm SOL}^{S}(K_{x_0}, f)$. By \cite[Theorem 4.1]{Kim3}, we directly obtain the assertions $(\rm ii)-(\rm iv)$. Conversely, if one of the conditions $(\rm ii)-(\rm iv)$ holds, then ${\rm SOL}^{S}(K_{x_0}, f)\neq\emptyset$ by \cite[Theorem 4.1]{Kim3} and $[f(K)]_{f(x_0)}=[f(K_{x_0})]_{f(x_0)}$. It follows from Lemma \ref{2023801} that ${\rm SOL}^{S}(K_{x_0}, f)\subseteq{\rm SOL}^{S}(K, f)$.  Hence, ${\rm SOL}^{S}(K, f)\neq\emptyset$. This proof is completed.\cvd



By Corollary \ref{C1}  and the proof of Theorem \ref{gpro} and Theorem  \ref{20231227ex}, we have the following results.

\begin{corollary}\label{240506}
	The following assertions are equivalent:
    \begin{itemize}
	\item[\rm(i)] $\rm{VOP}$$(K,f)$ admits a weak Pareto efficient solution.
	\item[\rm(ii)] There exist a nonempty index set $I\subseteq\{1,2,\dots,s\}$ and $x_0\in K$ such that $f$ is  weakly section-bounded from below on $K$ with respect to $I$ and $x_0$, and the inclusions that appear in (ii), (iii) and (iv) of Theorems  \ref{gpro} and \ref{20231227ex}  hold.
    \end{itemize}
\end{corollary}
{\it Proof}\,\,

``$(\rm ii)\Rightarrow (\rm i)$": It is a direct consequence of Theorem \ref{gpro} and Theorem  \ref{20231227ex}.

``$(\rm i)\Rightarrow (\rm ii)$":
By  Corollary \ref{C1}, there exists $x_0\in K$ such that the index set $I:=\{i\in \{1,2,\dots,s \}| f_{i}(x)\equiv f_{i}(x_0),\forall x\in K_{x_0}\}$ is nonempty and $x_0\in {\rm SOL}^{S}(K_{x_0}, f_I)$. Clearly, $f$ is  weakly section-bounded from below on $K$ with respect to $I$ and $x_0$. By  the proof of  Theorem \ref{gpro}, we  know that the inclusions that appear in (ii), (iii) and (iv) of Theorem \ref{gpro} hold. Similarly, by the proof of  Theorem \ref{20231227ex}, we have that the inclusions that appear in (ii), (iii) and (iv) of Theorem  \ref{20231227ex} also hold.    \cvd

 In addition, by Corollary \ref{240506}, we can obtain the following result.
\begin{corollary}\label{240707}
	Assume that $\rm{VOP}$$(K,f)$ admits a weak Pareto efficient solution. Then there exist a nonempty index set $I\subseteq\{1,2,\dots,s\}$ and $x_0\in K$ such that  the inclusions that appear in (ii), (iii) and (iv) of Theorem  \ref{20231227ex}   hold if and only if the inclusions  that   appear in (ii), (iii) and (iv) of Theorem \ref{gpro} hold, respectively.
\end{corollary}

\begin{remark}\label{re24707}    In order to find the point $x_0$ and the set $I$ which appear in  Corollary \ref{240506} and Corollary \ref{240707},
we only consider $x_0\in SOL^W(K, f)$  and $I\subseteq \{i\in \{1,2,\dots,s \}| f_{i}(x)\equiv f_{i}(x_0),\forall x\in K_{x_0}\}\ne\emptyset$.
 A natural question arises when the problem $\rm{VOP}$$(K,f)$ admits a weak Pareto efficient solution:   for any nonempty index set $I\subseteq\{1,2,\dots,s\}$ and any $x_0\in K$ such that $f$ is  weakly section-bounded from below on $K$ with respect to $I$ and $x_0$, do the respective inclusions that appear in (ii), (iii) and (iv) of Theorems \ref{gpro} and \ref{20231227ex} hold at the same time?
The following example provides a negative answer to the above question showing that when $SOL^W(K, f)\neq\emptyset$  and $f$ is  weakly section-bounded from below on $K$ with respect to some $I$ and $x_0$,  the inclusions  (ii) and (iii) of Theorem \ref{gpro} may hold, but the respective inclusions   (ii) and (iii) of Theorem \ref{20231227ex} may  not hold.
\end{remark}

\begin{example}\label{6.29}
	Let $f:=(f_1, f_2): \mathbf{R}^2\to \mathbf{R}^2$ with $f_1(x)=\frac{1}{2} x_1^2x_2+x_1$, $f_2(x)=|x_1|$, $K:=\{x\in\mathbf{R}^2\vert x_1\geq -1, x_2\geq 0\}$ and $x_0=(x_1, x_2)=(-1, 2)$.
Then $K_{x_0}:=\{x\in K\vert f_1(x)=\frac{1}{2} x_1^2x_2+x_1\leq 0, f_2(x)=|x_1|\leq 1\}$. Clearly, the set $K_{x_0}$ is closed, $f:=(f_{1}, f_{2})$ is  a locally Lipschitz mapping, $f$ is  weakly section-bounded from below on $K$ with respect to $I=\{1\}$ and $x_0=(-1, 2)$ (since $f_I=f_1\in [-1, 0]$ on $K_{x_0}$). Let $x_n=(-\frac{1}{n}, n)$ with $n\geq 1$.  We observe that $\{x_n\}\subseteq K_{x_0}\subseteq K$, $f(x_n)\leq f(x_0)$ with $n\geq 1$ and $\|x_n\|\to +\infty$ as $n\to +\infty$. Since the points $x_n$ with $n\geq 1$ are all inner points of $K$ and $K_{x_0}$,  we can easily  calculate $\nu(x_n)\leq\frac{1}{2n^2}$ and $\nu_{I, K_{x_0}}(x_n)=\frac{1}{2n^2}$. Therefore, $\nu(x_n)\to 0$ and  $\nu_{I, K_{x_0}}(x_n)\to 0$ as $n\to +\infty$, yielding $\|x_n\|\nu(x_n)\to 0$ and  $\|x_n\|\nu_{I, K_{x_0}}(x_n)\to 0$ as $n\to +\infty$. Because $f_I(x_n)=f_1(x_n)=-\frac{1}{2n}\to 0$ as $n\to +\infty$, by definitions,  we have $0\in \widetilde{K}^I_{\infty, \leq f(x_0)}(f, K)$, $0\in \widetilde{K}_{\infty, \leq f_I(x_0)}(f_I, K_{x_0})$, $0\in K^{I}_{\infty,\leq f(x_0)}(f, K)$ and $0\in K_{\infty,\leq f_I(x_0)}(f_I, K_{x_0})$. We can easily see that, for any  sequence $\{y_n=(y^1_n, y^2_n)\}\subseteq K_{x_0}$ with $\|y_n\|\to +\infty$ as $n\to +\infty$, since $f_I(y_n)=f_1(y_n)\in [-1, 0]$ for each $n$, then  $y^1_n\to 0$ as $n\to +\infty$. By $y^1_n\leq\frac{1}{2} (y^1_n)^2y^2_n+y^1_n\leq 0$, we have $f_I(y_n)=f_1(y_n)\to 0$ as $n\to +\infty$. Thus, we  obtain that $\widetilde{K}^I_{\infty, \leq f(x_0)}(f, K)=\widetilde{K}_{\infty, \leq f_I(x_0)}(f_I, K_{x_0})=\{0\}$ and $K^{I}_{\infty,\leq f(x_0)}(f, K)= K_{\infty,\leq f_I(x_0)}(f_I, K_{x_0})=\{0\}$. It is known that $\partial f_2((0, 0))=\{v=(v_1, 0)\in \mathbf{R}^2\vert v_1\in [-1,1]\}$(by Lemma \ref{classical subdifferential} and  \cite[Proposition 10.5]{Rock2}). Thus, we can see that there exist $v\in \partial f_2((0, 0))$ and $\alpha=(\alpha_1, \alpha_2)\in \mathbf{R}^2_+$ with $\sum^{2}_{i=1}\alpha_i=1$ such that $\alpha_1\nabla f_1(x_0)+\alpha_2v=\alpha_1(1, 0)+\alpha_2v=0$, yielding $\nu((0, 0))=0$. It follows from $(0, 0)\in K_{x_0}$ and $f_I((0, 0))=f_1((0, 0))=0$ that $0\in K^I_{0, \leq f(x_0)}(f, K)$. Therefore, we have $\widetilde{K}^I_{\infty, \leq f(x_0)}(f, K)\subseteq K^I_{0, \leq f(x_0)}(f, K)$ and $K^{I}_{\infty,\leq f(x_0)}(f, K)\subseteq K^I_{0, \leq f(x_0)}(f, K)$: by Theorem \ref{gpro} we have $SOL^W(K, f)\neq\emptyset$.  However, we can easily prove that $\nu_{I,K_{x_0}}(x)\neq 0$ for any $x\in K_{x_0}$, so that we have $0\notin K_{0, \leq f_I(x_0)}(f_I, K_{x_0})$. It follows that $\widetilde{K}_{\infty, \leq f_I(x_0)}(f_I, K_{x_0})\not\subseteq K_{0, \leq f_I(x_0)}(f_I, K_{x_0})$ and $K_{\infty,\leq f_I(x_0)}(f_I, K_{x_0})\not\subseteq K_{0, \leq f_I(x_0)}(f_I, K_{x_0})$.

\end{example}

The following example shows that when $SOL^W(K, f)\neq\emptyset$  and $f$ is  weakly section-bounded from below on $K$ with respect to $I$ and $x_0$, we have that the inclusion in (iv) of Theorem \ref{20231227ex} holds, but the inclusion in (iv) of Theorem \ref{gpro} does not hold.

\begin{example}\label{ex247071}
	Let $f:=(f_1, f_2): \mathbf{R}^2\to \mathbf{R}^2$ with $f_1(x)=(x_1-x_2)^2+x_1-x_2$, $f_2(x)=-x_1$, $K:=\{x\in\mathbf{R}^2\vert x_1\geq 0, x_2\geq 0\}$ and $x_0=(x_1, x_2)=(0, 0)$.
Then $K_{x_0}:=\{x\in K\vert (x_1-x_2)^2+x_1-x_2\leq 0\}$. Clearly, the set $K_{x_0}$ is closed and $f:=(f_{1}, f_{2})$ is  a locally Lipschitz mapping and $f$ is  weakly section-bounded from below on $K$ with respect to $I=\{1\}$ and $x_0=(0, 0)$ (since $f_1\in [-\frac{1}{2}, 0]$ on $K_{x_0}$).
It is not difficult to  prove that $T_{\infty, \leq f_I(x_0)}(f_I, K_{x_0})=\{-\frac{1}{2}\}\subseteq K_{0, \leq f_I(x_0)}(f_I, K_{x_0})$. Thus, we have $SOL^W(K, f)\neq\emptyset$ by Theorem \ref{20231227ex}. Let $x_n=(n, n)$ with $n\geq 1$. Then $\{x_n\}\subseteq K_{x_0}\subseteq K$, $f(x_n)\leq f(x_0)$ with $n\geq 1$ and $\|x_n\|\to +\infty$ as $n\to +\infty$. It can be  shown that the sequence  $\{x_n\}\subseteq\Gamma(f, K)$. By $\{x_n\}\subseteq K_{x_0}$ and $f_I(x_n)=f_1(x_n)\equiv 0$, we have  $0\in T^I_{\infty, \leq f(x_0)}(f, K)$. However, we can easily prove that $0\notin K^I_{0, \leq f(x_0)}(f, K)$, so that we have $T^I_{\infty, \leq f(x_0)}(f, K)\not\subseteq K^I_{0, \leq f(x_0)}(f, K)$.

Finally, it is worth noting that both Examples \ref{6.29} and  \ref{ex247071} also show that under the assumption that there exist a nonempty index set $I\subseteq\{1,2,\dots,s\}$ and $x_0\in K$ such that $f$ is  weakly section-bounded from below on $K$ with respect to $I$ and $x_0$, the inclusions that appear in (ii), (iii) and (iv) of Theorem  \ref{20231227ex} may not be equivalent to the inclusions  that  appear in (ii), (iii) and (iv) of Theorem \ref{gpro}, respectively.

\end{example}

\begin{remark}
As a final simple observation we may mention that   statements
(ii), (iii) and (iv) of Theorem  \ref{20231227ex} are equivalent to  (ii), (iii) and (iv) of Theorem \ref{gpro}, while
 the set $I$ and the point $x_0\in K$ are independently chosen in each statement.
\end{remark}

\section{Conclusion}
  As we know, a characteristic of the nonemptiness of Pareto efficient solutions set for a constrained vector optimization problem is that the image of the objective mapping has a bounded section and this result is easy to prove (see, e.g., Proposition \ref{I=123s}). However, from this viewpoint, it may be difficult to obtain a similar property related to the nonemptiness of weak Pareto efficient solution set (see Proposition \ref{property of nonempty}). In this paper, firstly, a new characterization of nonemptiness of the weak Pareto efficient solution set of  the constrained vector optimization problem is obtained. Secondly, based on the property of weak section-boundedness from below, we discuss the relationships between the notions of properness, Palais-Smale, weak Palais-Smale and M-tameness conditions with respect to some index set $I$ for the restricted vector mapping on the constraint set. Finally, by using the new characterization of nonemptiness of weak Pareto efficient solution set and a novel approach based on variational analysis and generalized differentiation, we derive some new  sufficient and necessary conditions for the existence of weak Pareto efficient solutions of nonsmooth constrained problems with locally Lipschitz mappings. The obtained results improve  Theorem 3.1 in \cite{Kim3} and develop some recent existence results of Pareto efficient solutions (e.g., \cite [Theorem 4.1]{Kim3}).



\begin{thebibliography}{10}

\bibitem{HHV} Ha, H.V., Pham, T.S.(eds.): Genericity in Polynomial Optimization. World Science Publishing, Singapore (2017)

\bibitem{BR} Benedetti, R., Risler, J.J.(eds.): Real Algebraic and Semi-Algebraic Sets. Hermann, Paris (1990)

\bibitem{Kim2} Kim, D.S., Pham, T.S., Tuyen, N.V.: On the existence of Pareto solutions for polynomial vector optimization problems, Math. Program.  Ser. A, \textbf{177}, 321--341 (2019)

\bibitem{LGJ}  Lee, J.H., Sisarat, N., Jiao, L.G.: Multi-objective convex polynomial optimization and semidefinite programming relaxations, J. Global Optim., \textbf{80}, 117--138 (2021)

\bibitem{Duan} Duan, Y.R., Jiao, L.G., Wu, P.C., Zhou, Y.Y.: Existence of Pareto solutions for vector polynomial optimization problems with constraints, J. Optim. Theory Appl., \textbf{195}, 148--171 (2022)

\bibitem{LDY2} Liu, D.Y., Huang, L., Hu, R.: Frank-Wolfe type theorems for polynomial vector optimization problems, Pac. J. Optim., \textbf{19}, 463--475 (2023)

\bibitem{LDY} Liu, D.Y., Hu, R., Fang, Y.P.: Solvability of a regular polynomial vector optimization problem without convexity, Optim., \textbf{72}, 821--841 (2023)

\bibitem{Kim3} Kim, D.S., Mordukhovich, B.S., Pham, T.S., Tuyen, N.V.: Existence of efficient and properly efficient solutions to problems of constrained vector optimization, Math. Program. Ser. A, \textbf{190}, 259--283 (2021)

\bibitem{Kim5} Kim, D.S., Pham, T.S., Tuyen, N.V.: On the existence of Pareto solutions for polynomial vector optimization problems. Math. Program. Ser. A, \textbf{177}, 321--341 (2019)

\bibitem{BT2} Bao, B.T., Mordukhovich, B.S.: Relative Pareto minimizers for multiobjective problems: existence and optimality conditions, Math. Program. Ser. A, \textbf{122}, 301--347 (2010)

\bibitem{BT1} Bao, B.T., Mordukhovich, B.S.: Variational principles for set-valued mappings with applications to multiobjective optimization, Control Cybern., \textbf{36}, 531--562 (2007)

\bibitem{MK2} Mordukhovich, B.S.: Variational Analysis and Applications. Springer, Cham (2018)

\bibitem{Flores2} Flores-Baz\'{a}n, F., Vera, C.: On the set of weakly efficient minimizers for convex multiobjective programming, Oper. Res. Lett., \textbf{36}, 651--655 (2008)

\bibitem{Flores} Flores-Baz\'{a}n, F.: Ideal, weakly efficient solutions for vector optimization problems. Math. Program. Ser. A,  \textbf{93}, 453--475 (2002)

\bibitem{Guti} Guti\'{e}rrez, C., L\`{o}pez, R.: On the existence of weak efficient solutions of nonconvex vector optimization problems, J. Optim. Theory Appl., \textbf{185}, 880--902 (2020)

\bibitem{Jahn} Jahn, J.: Vector Optimization: Theory, Applications, and Extensions. Springer, Berlin (2004)

\bibitem{Luc} Luc, D.T.: Theory of Vector Optimization. Lecture Notes in Economics and Mathematical Systems. Springer, Berlin (1989)

\bibitem{Sawaragi} Sawaragi, Y., Nakayama, N., Tanino, T.: Theory of Multiobjective Optimization. Academic Press, New York (1985)

\bibitem{Borwein} Borwein, J.M.: On the existence of Pareto efficient points. Math. Oper. Res., \textbf{8}, 64--73 (1983)

\bibitem{Deng1} Deng, S.: Characterizations of the nonemptiness and compactness of solution sets in convex vector optimization, J. Optim. Theory Appl., \textbf{96}, 123--131 (1998)

\bibitem{Deng2} Deng, S.: Boundedness and nonemptiness of the efficient solution sets in multiobjective optimization, J. Optim. Theory Appl., \textbf{144}, 29--42 (2010)

\bibitem{Huy} Huy, N.Q., Kim, D.S., Tuyen, N.V.:  Existence theorems in vector optimization with generalized order, J. Optim. Theory Appl., \textbf{174}, 728--745 (2017)

\bibitem{Lara} Lara, F.: Generalized asymptotic functions in nonconvex multiobjective optimization problems. Optim., \textbf{66}, 1259--1272 (2017)

\bibitem{Luc2} Luc, D.T.: An existence theorem in vector optimization. Math. Oper. Res., \textbf{14}, 693--699 (1989)

\bibitem{0810ha} H\`a, T.X.D.: Variants of the Ekeland variational principle for a set-valued map involving the Clarke normal cone, J. Math. Anal. Appl., \textbf{316}, 346--356 (2006)

\bibitem{MK1} Mordukhovich, B.S.: Variational Analysis and Generalized Differentiation. Vol.1: Basic Theory, Vol.2: Applications. Springer, Berlin (2006)

\bibitem{Rock2} Rockafellar, R.T., Wets, R.J.-B.: Variational Analysis. Springer, Berlin (1998)

\bibitem{Rock1} Rockafellar, R.T.: Convex Analysis. Princeton University Press, Princeton  (1970)

\bibitem{Ekland} Ekeland, I.: Nonconvex minimization problems, Bull. Am. Math. Soc., \textbf{1}, 443--474 (1979)

\bibitem{CLark05} Clarke, F.H.: A new approach to Lagrange multipliers, Math. Oper. Res., \textbf{1}, 165--174 (1976)

\bibitem{Tammer} Tammer, Chr.: A generalization of Ekellandz's variational principle, Optim., \textbf{25}, 129--141 (1992)

\bibitem{Loridan} Loridan, P.: Necessary conditions for $\varepsilon$-optimality. In: Optimality and Stability in Mathematical Programming, Guignard, M.(ed), Mathematical Programming Studies \textbf{19}, 140--152, Springer Berlin, Heidelberg (1982)

\end{thebibliography}
\end{document}